\documentclass[11pt]{amsart}
\usepackage{amsmath,amssymb,amsthm}
\usepackage{hyperref}
\hypersetup{hidelinks=true}

\usepackage{bbm}

\usepackage[utf8]{inputenc}

\usepackage{graphicx}
\usepackage{enumerate}
\usepackage{amsfonts,amssymb,latexsym}
\usepackage{accents}
\usepackage{afterpage}
\usepackage[style=alphabetic,maxbibnames=99, backend=bibtex]{biblatex}
\usepackage{tikz}

\newcommand{\N}{\mathbb{N}}

\newcommand{\R}{\mathbb{R}}

\newcommand{\Z}{\mathbb{Z}}

\newcommand{\s}{\sum}

\newcommand\restr[2]{{
  \left.\kern-\nulldelimiterspace 
  #1 
  \vphantom{\big|} 
  \right|_{#2} 
  }}

\newcommand{\cco}{\overline{\operatorname{co}}}

\newtheorem{theorem}{Theorem}[section]
\newtheorem{lemma}[theorem]{Lemma}
\newtheorem{claim}[theorem]{Claim}
\newtheorem{prop}[theorem]{Proposition}
\newtheorem{corollary}[theorem]{Corollary}
\theoremstyle{definition}
\newtheorem{definition}[theorem]{Definition}

\theoremstyle{remark}
\newtheorem{remark}[theorem]{Remark}
\numberwithin{equation}{section}

\def\fnote#1{\footnote}

\def\R{{\mathbb R}}

\def\ignora#1{}
\def\n3#1{\left\vert  \! \left\vert \! \left\vert \, #1 \, \right\vert \!
  \right\vert \! \right\vert }

\include{draft}

\begin{document}



\title[Bases in nets]{Schauder bases in Lipschitz free spaces over nets in Banach spaces}

\author{Petr H\'ajek}\thanks{This research was supported by CAAS CZ.02.1.01/0.0/0.0/16-019/0000778
 and by the project  SGS21/056/OHK3/1T/13.}
\address[P. H\'ajek]{Czech Technical University in Prague, Faculty of Electrical Engineering.
Department of Mathematics, Technická 2, 166 27 Praha 6 (Czech Republic)}
\email{ hajek@math.cas.cz}

\author{ Rub\'en Medina}\thanks{The second author research has also been supported by MICINN (Spain) Project PGC2018-093794-B-I00 and MIU (Spain) FPU19/04085 Grant.}
\address[R. Medina]{Universidad de Granada, Facultad de Ciencias.
Departamento de An\'{a}lisis Matem\'{a}tico, 18071-Granada
(Spain); and Czech technical University in Prague, Faculty of Electrical Engineering.
Department of Mathematics, Technická 2, 166 27 Praha 6 (Czech Republic)}
\email{rubenmedina@ugr.es}
\urladdr{\url{https://www.ugr.es/personal/ae3750ed9865e58ab7ad9e11e37f72f4}}

\maketitle

\begin{abstract}

In the present note we give two explicit constructions (based on a retractional argument) of a Schauder basis for the Lipschitz free space $\mathcal{F}(N)$, 
over certain uniformly discrete metric spaces $N$.
The first one applies to every net $N$ in
a finite dimensional Banach space, leading to the basis constant independent of the dimension. The second one applies to grids in Banach spaces with an FDD.

As a corollary, we obtain a retractional Schauder basis for the Lipschitz free space $\mathcal{F}(N)$ over a net $N$ in every Banach space $X$ with a Schauder basis containing a copy of $c_0$, as well as
in every Banach space with a $c_0$-like FDD.\end{abstract}

\section{Introduction}

The structural properties of metric spaces, in terms of the various extensional or retractional properties of their subsets and Lipschitz functions thereon, are a very active field of research with many
applications in surrounding areas of Banach space geometry, theoretical computer science, geometric group theory and so on. A convenient conceptual framework in this respect is
provided by the notion of a Lipschitz free space $\mathcal{F}(M)$ over a metric space $M$ (also known as the Arens-Eells space). It is a Banach space containing an isometric copy of the
original metric space $M$, which allows the natural linearization of Lipschitz mappings between metric spaces. The structural properties of  Lipschitz free spaces $\mathcal{F}(M)$
are still very far from being understood, except for some very special cases of $M$ such as e.g.  the doubling metric spaces, or (separable) uniformly discrete spaces $M$.
In the former case it was shown by Lancien and Perneck\'a  \cite{LP13}  that $\mathcal{F}(M)$ has a $\lambda$-BAP for any $M\subset\R^n$, where $\lambda$ is of the order $\sqrt{n}$ and independent of $M$. It remains an 
important open problem whether  $\lambda$ can be chosen a constant independent of $n$ (e.g. \cite{PS15}).

 In the latter case,
Kalton \cite{Kal12} has shown that $\mathcal{F}(M)$ is a Schur space with the RNP property and  the approximation property (AP for short). Kalton  has posed an important
question whether these spaces have the bounded approximation property (BAP) as well. In fact, in order to solve this problem it suffices to consider the special case when $M$
is a net in a  Banach space $X$. 
 A positive answer in the separable case would imply that every separable Banach space is approximable in the sense defined in \cite{Kal12} while a counterexample would yield a renorming of $\ell_1$ without the 
metric approximation property (MAP), answering a classical open question.

Kalton went on proving that for a Banach space $X$ with a separable dual (or a separable dual space itself) the answer is positive. So for every net $N$ in such $X$
the Lipschitz free  space $\mathcal{F}(N)$ has the BAP.
 The  case of a general separable Banach space $X$ still remains open.
Let us add the remark that nets in metric spaces are quite useful (e.g. in theoretical computer science) and their study is an ongoing topic in the last years (see \cite{HN17}, \cite{DOS+08}, \cite{Kal12}).

Our note contains two main results which considerably strengthen the above mentioned results in some special cases. 

We construct an explicit family of Lipschitz retractions for nets in finite dimensional Banach spaces,
which yields a Schauder basis for their Lipschitz free spaces whose constant is independent of the dimension.
Based on this ingredient, we proceed with the construction of an explicit family of Lipschitz retractions (again leading to a Schauder basis for the free space) for  grids in Banach spaces with an FDD.
Our construction can be applied to obtain a retractional Schauder basis for  $\mathcal{F}(N)$, where $N$ is a net in several types of Banach spaces $X$.
In particular, $X$ can be chosen a separable Banach space with a Schauder basis, which contains a copy of $c_0$, or a 
 Banach space admitting a $c_0$-like FDD.
Let us now proceed with the necessary backgroud.

Given $a,b\in\R^+$ we will say that a subset $N$ of a metric space $M$ is an $(a,b)$-net of $M$ whenever it is $a$-separated and $b$-dense. More precisely, $N$ is considered to be $a$-separated whenever $d(x,y)\ge a$ for every $x,y\in N$, $x\neq y$, and $b$-dense if for every $x\in M$ there is $y\in N$ such that $d(x,y)\le b$. We will say that $N\subset M$ is a net whenever there are $a,b\in \R^+$ such that $N$ is an $(a,b)$-net.

A map $T$ from a metric space $M$ into another metric space $N$ is said to be Lipschitz if there exists some $\lambda>0$ such that
$$d(T(x),T(y))\le \lambda d(x,y) \;\;\;\;\forall x,y\in M.$$
We say that $\lambda$ is a Lipschitz constant for $T$ and we call the infimum of all Lipschitz constants for $T$ the Lipschitz norm of $T$, that is,
$$||T||_{\text{Lip}}=\sup\limits_{x,y\in M, x\ne y}\frac{d\big(T(x),T(y)\big)}{d(x,y)}.$$
If $\lambda>0$ is a Lipschitz constant for $T$ then we say that $T$ is $\lambda$-Lipschitz.

\begin{definition}
Let $N_1$ and $N_2$ be metric spaces. We say that $N_1$ and $N_2$ are Lipschitz equivalent whenever there is a bijection $F:N_1\to N_2$ such that both $F$ and $F^{-1}$ are Lipschitz.
We put $\text{dist}(F)=\|F\|_{\text{Lip}}\|F^{-1}\|_{\text{Lip}}$ the Lipschitz distorsion of $F$.
\end{definition}

\begin{figure}\label{fignet1}
\centering

\tikzset{every picture/.style={line width=0.75pt}} 

\begin{tikzpicture}[x=0.75pt,y=0.75pt,yscale=-0.8,xscale=0.8]

\draw  [fill={rgb, 255:red, 0; green, 0; blue, 0 }  ,fill opacity=1 ] (366.67,168.17) .. controls (366.67,166.79) and (367.79,165.67) .. (369.17,165.67) .. controls (370.55,165.67) and (371.67,166.79) .. (371.67,168.17) .. controls (371.67,169.55) and (370.55,170.67) .. (369.17,170.67) .. controls (367.79,170.67) and (366.67,169.55) .. (366.67,168.17) -- cycle ;
\draw  [fill={rgb, 255:red, 0; green, 0; blue, 0 }  ,fill opacity=1 ] (431.67,138.17) .. controls (431.67,136.79) and (432.79,135.67) .. (434.17,135.67) .. controls (435.55,135.67) and (436.67,136.79) .. (436.67,138.17) .. controls (436.67,139.55) and (435.55,140.67) .. (434.17,140.67) .. controls (432.79,140.67) and (431.67,139.55) .. (431.67,138.17) -- cycle ;
\draw  [fill={rgb, 255:red, 0; green, 0; blue, 0 }  ,fill opacity=1 ] (436.67,83.17) .. controls (436.67,81.79) and (437.79,80.67) .. (439.17,80.67) .. controls (440.55,80.67) and (441.67,81.79) .. (441.67,83.17) .. controls (441.67,84.55) and (440.55,85.67) .. (439.17,85.67) .. controls (437.79,85.67) and (436.67,84.55) .. (436.67,83.17) -- cycle ;
\draw  [fill={rgb, 255:red, 0; green, 0; blue, 0 }  ,fill opacity=1 ] (496.67,53.17) .. controls (496.67,51.79) and (497.79,50.67) .. (499.17,50.67) .. controls (500.55,50.67) and (501.67,51.79) .. (501.67,53.17) .. controls (501.67,54.55) and (500.55,55.67) .. (499.17,55.67) .. controls (497.79,55.67) and (496.67,54.55) .. (496.67,53.17) -- cycle ;
\draw  [fill={rgb, 255:red, 0; green, 0; blue, 0 }  ,fill opacity=1 ] (531.67,108.17) .. controls (531.67,106.79) and (532.79,105.67) .. (534.17,105.67) .. controls (535.55,105.67) and (536.67,106.79) .. (536.67,108.17) .. controls (536.67,109.55) and (535.55,110.67) .. (534.17,110.67) .. controls (532.79,110.67) and (531.67,109.55) .. (531.67,108.17) -- cycle ;
\draw  [fill={rgb, 255:red, 0; green, 0; blue, 0 }  ,fill opacity=1 ] (501.67,148.17) .. controls (501.67,146.79) and (502.79,145.67) .. (504.17,145.67) .. controls (505.55,145.67) and (506.67,146.79) .. (506.67,148.17) .. controls (506.67,149.55) and (505.55,150.67) .. (504.17,150.67) .. controls (502.79,150.67) and (501.67,149.55) .. (501.67,148.17) -- cycle ;
\draw  [fill={rgb, 255:red, 0; green, 0; blue, 0 }  ,fill opacity=1 ] (591.67,78.17) .. controls (591.67,76.79) and (592.79,75.67) .. (594.17,75.67) .. controls (595.55,75.67) and (596.67,76.79) .. (596.67,78.17) .. controls (596.67,79.55) and (595.55,80.67) .. (594.17,80.67) .. controls (592.79,80.67) and (591.67,79.55) .. (591.67,78.17) -- cycle ;
\draw  [fill={rgb, 255:red, 0; green, 0; blue, 0 }  ,fill opacity=1 ] (606.67,143.17) .. controls (606.67,141.79) and (607.79,140.67) .. (609.17,140.67) .. controls (610.55,140.67) and (611.67,141.79) .. (611.67,143.17) .. controls (611.67,144.55) and (610.55,145.67) .. (609.17,145.67) .. controls (607.79,145.67) and (606.67,144.55) .. (606.67,143.17) -- cycle ;
\draw  [fill={rgb, 255:red, 0; green, 0; blue, 0 }  ,fill opacity=1 ] (561.67,198.17) .. controls (561.67,196.79) and (562.79,195.67) .. (564.17,195.67) .. controls (565.55,195.67) and (566.67,196.79) .. (566.67,198.17) .. controls (566.67,199.55) and (565.55,200.67) .. (564.17,200.67) .. controls (562.79,200.67) and (561.67,199.55) .. (561.67,198.17) -- cycle ;
\draw  [fill={rgb, 255:red, 0; green, 0; blue, 0 }  ,fill opacity=1 ] (456.67,198.17) .. controls (456.67,196.79) and (457.79,195.67) .. (459.17,195.67) .. controls (460.55,195.67) and (461.67,196.79) .. (461.67,198.17) .. controls (461.67,199.55) and (460.55,200.67) .. (459.17,200.67) .. controls (457.79,200.67) and (456.67,199.55) .. (456.67,198.17) -- cycle ;
\draw    (375.67,164.93) -- (425.62,141.54) ;
\draw [shift={(428.33,140.27)}, rotate = 154.9] [fill={rgb, 255:red, 0; green, 0; blue, 0 }  ][line width=0.08]  [draw opacity=0] (10.72,-5.15) -- (0,0) -- (10.72,5.15) -- (7.12,0) -- cycle    ;
\draw    (438.33,90.27) -- (435.25,127.28) ;
\draw [shift={(435,130.27)}, rotate = 274.76] [fill={rgb, 255:red, 0; green, 0; blue, 0 }  ][line width=0.08]  [draw opacity=0] (10.72,-5.15) -- (0,0) -- (10.72,5.15) -- (7.12,0) -- cycle    ;
\draw    (503,58.93) -- (529.39,100.4) ;
\draw [shift={(531,102.93)}, rotate = 237.53] [fill={rgb, 255:red, 0; green, 0; blue, 0 }  ][line width=0.08]  [draw opacity=0] (10.72,-5.15) -- (0,0) -- (10.72,5.15) -- (7.12,0) -- cycle    ;
\draw    (441.67,139.6) -- (494.69,145.91) ;
\draw [shift={(497.67,146.27)}, rotate = 186.79] [fill={rgb, 255:red, 0; green, 0; blue, 0 }  ][line width=0.08]  [draw opacity=0] (10.72,-5.15) -- (0,0) -- (10.72,5.15) -- (7.12,0) -- cycle    ;
\draw    (530.33,114.27) -- (511.42,140.5) ;
\draw [shift={(509.67,142.93)}, rotate = 305.79] [fill={rgb, 255:red, 0; green, 0; blue, 0 }  ][line width=0.08]  [draw opacity=0] (10.72,-5.15) -- (0,0) -- (10.72,5.15) -- (7.12,0) -- cycle    ;
\draw    (464.33,192.93) -- (498.31,155.81) ;
\draw [shift={(500.33,153.6)}, rotate = 132.47] [fill={rgb, 255:red, 0; green, 0; blue, 0 }  ][line width=0.08]  [draw opacity=0] (10.72,-5.15) -- (0,0) -- (10.72,5.15) -- (7.12,0) -- cycle    ;
\draw    (587,81.6) -- (543.66,104.21) ;
\draw [shift={(541,105.6)}, rotate = 332.45] [fill={rgb, 255:red, 0; green, 0; blue, 0 }  ][line width=0.08]  [draw opacity=0] (10.72,-5.15) -- (0,0) -- (10.72,5.15) -- (7.12,0) -- cycle    ;
\draw    (559,194.27) -- (511.95,154.21) ;
\draw [shift={(509.67,152.27)}, rotate = 40.41] [fill={rgb, 255:red, 0; green, 0; blue, 0 }  ][line width=0.08]  [draw opacity=0] (10.72,-5.15) -- (0,0) -- (10.72,5.15) -- (7.12,0) -- cycle    ;
\draw    (603,140.53) -- (543.05,112.47) ;
\draw [shift={(540.33,111.2)}, rotate = 25.08] [fill={rgb, 255:red, 0; green, 0; blue, 0 }  ][line width=0.08]  [draw opacity=0] (10.72,-5.15) -- (0,0) -- (10.72,5.15) -- (7.12,0) -- cycle    ;
\draw  [fill={rgb, 255:red, 0; green, 0; blue, 0 }  ,fill opacity=1 ] (31.33,168.83) .. controls (31.33,167.45) and (32.45,166.33) .. (33.83,166.33) .. controls (35.21,166.33) and (36.33,167.45) .. (36.33,168.83) .. controls (36.33,170.21) and (35.21,171.33) .. (33.83,171.33) .. controls (32.45,171.33) and (31.33,170.21) .. (31.33,168.83) -- cycle ;
\draw  [fill={rgb, 255:red, 0; green, 0; blue, 0 }  ,fill opacity=1 ] (96.33,138.83) .. controls (96.33,137.45) and (97.45,136.33) .. (98.83,136.33) .. controls (100.21,136.33) and (101.33,137.45) .. (101.33,138.83) .. controls (101.33,140.21) and (100.21,141.33) .. (98.83,141.33) .. controls (97.45,141.33) and (96.33,140.21) .. (96.33,138.83) -- cycle ;
\draw  [fill={rgb, 255:red, 0; green, 0; blue, 0 }  ,fill opacity=1 ] (101.33,83.83) .. controls (101.33,82.45) and (102.45,81.33) .. (103.83,81.33) .. controls (105.21,81.33) and (106.33,82.45) .. (106.33,83.83) .. controls (106.33,85.21) and (105.21,86.33) .. (103.83,86.33) .. controls (102.45,86.33) and (101.33,85.21) .. (101.33,83.83) -- cycle ;
\draw  [fill={rgb, 255:red, 0; green, 0; blue, 0 }  ,fill opacity=1 ] (161.33,53.83) .. controls (161.33,52.45) and (162.45,51.33) .. (163.83,51.33) .. controls (165.21,51.33) and (166.33,52.45) .. (166.33,53.83) .. controls (166.33,55.21) and (165.21,56.33) .. (163.83,56.33) .. controls (162.45,56.33) and (161.33,55.21) .. (161.33,53.83) -- cycle ;
\draw  [fill={rgb, 255:red, 0; green, 0; blue, 0 }  ,fill opacity=1 ] (196.33,108.83) .. controls (196.33,107.45) and (197.45,106.33) .. (198.83,106.33) .. controls (200.21,106.33) and (201.33,107.45) .. (201.33,108.83) .. controls (201.33,110.21) and (200.21,111.33) .. (198.83,111.33) .. controls (197.45,111.33) and (196.33,110.21) .. (196.33,108.83) -- cycle ;
\draw  [fill={rgb, 255:red, 0; green, 0; blue, 0 }  ,fill opacity=1 ] (166.33,148.83) .. controls (166.33,147.45) and (167.45,146.33) .. (168.83,146.33) .. controls (170.21,146.33) and (171.33,147.45) .. (171.33,148.83) .. controls (171.33,150.21) and (170.21,151.33) .. (168.83,151.33) .. controls (167.45,151.33) and (166.33,150.21) .. (166.33,148.83) -- cycle ;
\draw  [fill={rgb, 255:red, 0; green, 0; blue, 0 }  ,fill opacity=1 ] (256.33,78.83) .. controls (256.33,77.45) and (257.45,76.33) .. (258.83,76.33) .. controls (260.21,76.33) and (261.33,77.45) .. (261.33,78.83) .. controls (261.33,80.21) and (260.21,81.33) .. (258.83,81.33) .. controls (257.45,81.33) and (256.33,80.21) .. (256.33,78.83) -- cycle ;
\draw  [fill={rgb, 255:red, 0; green, 0; blue, 0 }  ,fill opacity=1 ] (271.33,143.83) .. controls (271.33,142.45) and (272.45,141.33) .. (273.83,141.33) .. controls (275.21,141.33) and (276.33,142.45) .. (276.33,143.83) .. controls (276.33,145.21) and (275.21,146.33) .. (273.83,146.33) .. controls (272.45,146.33) and (271.33,145.21) .. (271.33,143.83) -- cycle ;
\draw  [fill={rgb, 255:red, 0; green, 0; blue, 0 }  ,fill opacity=1 ] (226.33,198.83) .. controls (226.33,197.45) and (227.45,196.33) .. (228.83,196.33) .. controls (230.21,196.33) and (231.33,197.45) .. (231.33,198.83) .. controls (231.33,200.21) and (230.21,201.33) .. (228.83,201.33) .. controls (227.45,201.33) and (226.33,200.21) .. (226.33,198.83) -- cycle ;
\draw  [fill={rgb, 255:red, 0; green, 0; blue, 0 }  ,fill opacity=1 ] (121.33,198.83) .. controls (121.33,197.45) and (122.45,196.33) .. (123.83,196.33) .. controls (125.21,196.33) and (126.33,197.45) .. (126.33,198.83) .. controls (126.33,200.21) and (125.21,201.33) .. (123.83,201.33) .. controls (122.45,201.33) and (121.33,200.21) .. (121.33,198.83) -- cycle ;
\draw   (27,40.53) -- (297.67,40.53) -- (297.67,210.53) -- (27,210.53) -- cycle ;
\draw   (362.33,40.53) -- (633,40.53) -- (633,210.53) -- (362.33,210.53) -- cycle ;
\draw  [fill={rgb, 255:red, 0; green, 0; blue, 0 }  ,fill opacity=1 ] (47.33,52.5) .. controls (47.33,51.12) and (48.45,50) .. (49.83,50) .. controls (51.21,50) and (52.33,51.12) .. (52.33,52.5) .. controls (52.33,53.88) and (51.21,55) .. (49.83,55) .. controls (48.45,55) and (47.33,53.88) .. (47.33,52.5) -- cycle ;
\draw  [fill={rgb, 255:red, 0; green, 0; blue, 0 }  ,fill opacity=1 ] (381.33,52.5) .. controls (381.33,51.12) and (382.45,50) .. (383.83,50) .. controls (385.21,50) and (386.33,51.12) .. (386.33,52.5) .. controls (386.33,53.88) and (385.21,55) .. (383.83,55) .. controls (382.45,55) and (381.33,53.88) .. (381.33,52.5) -- cycle ;
\draw    (389.67,56.53) -- (431.69,79.11) ;
\draw [shift={(434.33,80.53)}, rotate = 208.25] [fill={rgb, 255:red, 0; green, 0; blue, 0 }  ][line width=0.08]  [draw opacity=0] (10.72,-5.15) -- (0,0) -- (10.72,5.15) -- (7.12,0) -- cycle    ;
\draw [line width=0.75]    (308.02,122.37) -- (344.69,122.9)(307.98,125.37) -- (344.65,125.9) ;
\draw [shift={(353.67,124.53)}, rotate = 180.84] [fill={rgb, 255:red, 0; green, 0; blue, 0 }  ][line width=0.08]  [draw opacity=0] (14.29,-6.86) -- (0,0) -- (14.29,6.86) -- cycle    ;

\draw (371.17,168.67) node [anchor=north west][inner sep=0.75pt]  [font=\small] [align=left] {$\displaystyle x_{10}$};
\draw (439.97,62.86) node [anchor=north west][inner sep=0.75pt]  [font=\small] [align=left] {$\displaystyle x_{9}$};
\draw (441.97,119.2) node [anchor=north west][inner sep=0.75pt]  [font=\small] [align=left] {$\displaystyle x_{4}$};
\draw (447.97,174.2) node [anchor=north west][inner sep=0.75pt]  [font=\small] [align=left] {$\displaystyle x_{3}$};
\draw (371,58.33) node [anchor=north west][inner sep=0.75pt]  [font=\small] [align=left] {$\displaystyle x_{11}$};
\draw (493,123.4) node [anchor=north west][inner sep=0.75pt]  [font=\small]  {$x_{1}$};
\draw (532.33,79.67) node [anchor=north west][inner sep=0.75pt]  [font=\small]  {$x_{5}$};
\draw (509,45.4) node [anchor=north west][inner sep=0.75pt]  [font=\small]  {$x_{8}$};
\draw (599,57.4) node [anchor=north west][inner sep=0.75pt]  [font=\small]  {$x_{7}$};
\draw (602,118.4) node [anchor=north west][inner sep=0.75pt]  [font=\small]  {$x_{6}$};
\draw (562,174.4) node [anchor=north west][inner sep=0.75pt]  [font=\small]  {$x_{2}$};

\end{tikzpicture}
\caption{Giving a retractional basis to a net}
\end{figure}
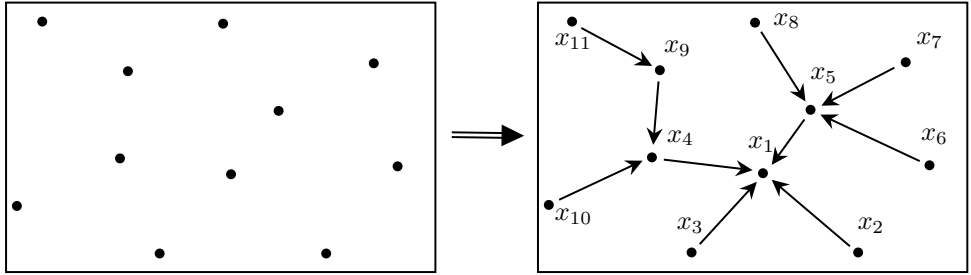

In this note, we are mainly interested in finding a retractional structure in nets. To this end, we are going to order the elements of the net
into a sequence and introduce $K$-Lipschitz retractions onto the first $n$ points of the sequence with a conmuting behaviour, as shown in figure \ref{fignet1}. More precisely,

\begin{definition}
Let $M$ be a countable metric space and $K\ge1$, a $K$-retractional basis of $M$ is a sequence of retractions $\{\varphi_n\}_{n=1}^\infty$, $\varphi_n:M\to M$ satisfying,
\begin{enumerate}
\item $\varphi_n(M)=M_n:=\bigcup_{j=1}^{n}\{x_j\}$,
\item $\bigcup_{j=1}^\infty\{x_j\}=M$,
\item $\varphi_n$ is $K$-Lipschitz for every $n\in\N$,
\item $\varphi_m\circ\varphi_n=\varphi_{\min(m,n)}$ for every $m,n\in\N$.
\end{enumerate}
We will say that $M$ has a retractional basis whenever it has a $K$-retractional basis for some $K\ge1$.
\end{definition}

It is clear that these retractions can be used to extend Lipschitz maps from a finite number of points $\{x_1,\dots,x_n\}$ to the whole net by composition, keeping the Lipschitz constant under control independently of the number of points from which the extension is made. By contrast, as shown by Naor and Rabani  \cite{NR17} for every $n\in\N$ there is an $n$-point subset $S$ of a metric space $M$ and a 1-Lipschitz function $f$ from $S$ into a certain Banach space $Z$ such that the Lipschitz constant of every extension of $f$ to the domain $M$ is bounded below by a multiple of $\sqrt{log\,n}$. 

The existence of retractional bases has been studied previously (see \cite{HN17} and \cite{Nov20}). Let us recall that by \cite{BK98} nets of a finite dimenisional space are not necessarily Lipschitz equivalent, so the problem is non-trivial
even in the case of a fixed normed space (and varying nets).

It is easy to see that retractional bases in nets are clearly preserved under Lipschitz equivalences. More precisely,

\begin{prop}\label{lipequiv}
Let $N_1$ and $N_2$ be countable Lipschitz equivalent metric spaces, with distorsion $D$. If $N_1$ has a $K$-retractional basis for some $K>0$ then $N_2$ has a $DK$-retractional basis. 
\end{prop}

A Schauder basis for a real Banach space $X$ is a sequence $(x_n)\subset X$ with the property that for every $x\in X$, there exists a unique sequence $(\alpha_n)\subset \R$ such that
$$\bigg|\bigg| x-\sum\limits_{i=1}^n\alpha_ix_i \bigg|\bigg|\xrightarrow{n\to\infty}0.$$
If $(x_n)$ is a Schauder basis then there is a $K>0$ such that $\Big|\Big|\sum\limits_{i=1}^n\alpha_ix_i \Big|\Big|\le K||x||$ for every $x\in X$ and $n\in\N$. In this case we say that $(x_n)$ is a $K$-Schauder basis.
A more general concept, also used in our note is the following.

\begin{definition}
A sequence $(X_n)$ of finite dimensional subspaces of a Banach space $X$ is called the finite dimensional decomposition (FDD for short) if for every $x\in X$ there is a unique sequence $x_n\in X_n$ so that
$$\bigg|\bigg| x-\sum\limits_{i=1}^nx_i \bigg|\bigg|\xrightarrow{n\to\infty}0.$$
In this case we say that $X$ has an FDD and we call the projections $P_n(x)=\sum\limits_{i=1}^nx_i$, which are again uniformly bounded, the natural projections of the FDD.
\end{definition}

A retractional basis for a net implies that the Lipschitz-free space over the net has a Schauder basis, made out of the linearization of the retractions.

\begin{theorem}[\cite{HN17}]\label{theoretr}
Let $M$ be a countable metric space with a $K$-retractional basis  $\{\varphi_n\}_{n=1}^\infty$. Then the Lipschitz-free space over $M$ has a $K$-Schauder basis.
\end{theorem}

In \cite{HN17}  a retractional basis in grids of spaces with an unconditional basis (Theorem 13) was constructed.  In section \ref{FDD} we  generalize these results to spaces with a Schauder basis using a completely different approach.

The paper is divided into two main Sections. In Section  \ref{finite} we tackle the problem in finite dimensional spaces. We first define a family of nets in $\R^n$ (we call those nets spiderwebs) such that every net of $\R^n$ is Lipschitz equivalent to a spider web with a distorsion independent of $n$. Then, we construct a $K$-retractional basis for an arbitrary $(a,b)$-net of the family with $K$ depending only on $\frac ba$, leading to our first main result Theorem \ref{theofinite}. In Section \ref{FDD}  we focus on grids of spaces with a Schauder decomposition (FDD). We prove our  main technical result Theorem \ref{theoFDD} showing that if $X$ is a space with an FDD, then some grids of the FDD have a retractional basis (a precise definition of this kind of grids is given in Section \ref{FDD}). 
 This result, combined with the coefficient quantization technique in \cite{DOS+08} yields our second main result, that every net of a Banach space with a Schauder basis containing a copy of $c_0$,
or a Banach space with a $c_0$-like FDD has a retractional basis (Corollary \ref{corFDD}). 

Let us remark that the importance of the mentioned coefficient quantization approach can be gleaned from 
 \cite{DOS+08}, \cite{CDO+08} and \cite{DD03}. In particular, in \cite{DD03} Fourier series of bandlimited functions (functions whose Fourier transform has a bounded support) are approximated replacing Fourier series coefficients by quantized coefficients, obtaining a net using a two-sided version of the $\Sigma-\Delta$ quantization algorithm. It is proved in the Main Theorem of \cite{DOS+08} that there is a Schauder basis of a space such that every grid with respect to that basis is a net if and only if the space has a Schauder basis and contains a copy of  $c_0$.

For background on the approximation properties of Banach spaces we refer to the Handbook article \cite{Cas01}. We will use the standard notation and terminology of the Banach spaces theory, as in \cite{Fab1}.
For background on Lipschitz-free spaces and its main properties we refer to \cite{GK03}.

\section{Finite dimensional nets}\label{finite}

Let $X$ be a finite dimensional Banach space. Our aim in this section is to prove that every $(a,b)$-net of $X$ has a retractional  basis with a bound depending only on
the fraction $\frac ba$ (and independent of the dimension of $X$). To this end, we will reduce the problem to the special case of nets described below.

\begin{definition}
Let $a,b\in\R^+$. We say that $N\subset X$ is a base of an $(a,b)$-spiderweb of $X$ whenever there exists a sequence of nonempty subsets $S_n\subset nS_X$ with $n\in\N\cup\{0\}$ such that
$$\begin{cases}2S_{2^n}\subset S_{2^{n+1}}\;\;\forall n \in\N\cup\{0\},\\S_n\text{ is }(a,b)\text{-net of }nS_X\;\;\;\;\forall n\in\N,\\
S_m=\frac{m}{2^n}S_{2^n}\;\;\forall m\in\{2^n+1,\dots,2^{n+1}-1\}\;\;\forall n\in\N,\end{cases}$$
 where $N=\bigcup\limits_{n=0}^{\infty}S_n$.
\end{definition}

The existence in $X$ of a base of an $(a,2a)$-spiderweb  for every $a\in (0,2)$ is proved by induction, taking $S_1$ an $(a,a)$-net of $S_X$ and finding $S_{2^n}$ as a maximal $a$-separated subset of $2^nS_X$ containing $2S_{2^{n-1}}$.

\begin{definition}
Let $a,b\in\R^+$. We say that $S\subset X$ is an $(a,b)$-spiderweb of $X$ if it is an $(a,b)$-net  containing as a subset a base of an $(a,b)$-spiderweb of $X$.
\end{definition}

\begin{lemma}\label{spider}
Every $(a,b)$-net is Lipschitz equivalent to an $\big(\frac{a}{6b},2\big)$-spiderweb with a distorsion at most $\big(\frac{2b}{a}+1\big)\big(\frac{4b}{a}+1\big)$.
\begin{proof}
Clearly, every  $(a, b)$-net is bi-Lipschitz equivalent, with a distorsion $1$, to an
 $(\widetilde a,\widetilde b)$-net with $\widetilde b=1/3$ and $\widetilde a=\frac{a}{3b}$.
Let $\widetilde{N}\subset X$ be a base of a $(1,2)$-spiderweb for $X$ and let $N\subset X$ be an $(\widetilde a,\widetilde b)$-net for such $\widetilde a$ and $\widetilde b$. We consider $T:\widetilde{N}\to N$ the nearest point map
(which is not uniquely determined). First, let us see that $T$ is injective. By contradiction, if $x,y\in \widetilde{N}$ are distinct and such that $T(x)=T(y)=z\in N$ then $||x-y||\le||x-z||+||y-z||\le 2/3$ which is impossible as $\widetilde{N}$ is $1$-separated. We define the set $S=\widetilde{N}\cup(N\setminus T(\widetilde{N}))$ and the map $F:N\to S$  as
$$F(x)=\begin{cases} x\;\;\;&\text{if }x\in N\setminus T(\widetilde{N}),\\ T^{-1}(x)&\text{if }x\in T(\widetilde{N}). \end{cases}$$
$F$ is well-defined and it is a bijection satisfying $||F(x)-x||\le \widetilde b=1/3$ for every $x\in N$. We are going to prove that $S$ is an $(\widetilde a /2,2)$-spiderweb. Obviously, $S$ contains $\widetilde{N}$ so it  remains to prove that $S$ is an $(\widetilde a /2,2)$-net. Let us first prove that $S$ is $\widetilde a /2$-separated. Take $x,y\in S$. If $x,y\in N\setminus T(\widetilde{N})$ then obviously $||x-y||\ge \widetilde a>\widetilde a/2$ and if $x,y\in\widetilde{N}$ then clearly $||x-y||\ge3\widetilde b\ge \widetilde a/2$. It  remains to deal with the case when $x\in N\setminus T(\widetilde{N})$ and $y\in\widetilde{N}$. In this case, assuming by contradiction that $||x-y||<\widetilde a /2$, we have $||T(y)-x||\le||T(y)-y||+||y-x||<\widetilde a /2+\widetilde a /2=\widetilde a$ which is again impossible because $N$ is $\widetilde a$-separated. Finally, let us prove that $S$ is $2$-dense in $X$. In fact, it is $2/3$-dense since for every $x\in X$ there
 exists a $y\in S$ such that $||x-F^{-1}(y)||\le \widetilde b$. As $N$ is $\widetilde b$-dense  we have that $||x-y||\le ||x-F^{-1}(y)||+||F^{-1}(y)-y||\le \widetilde b+\widetilde b=2\widetilde b=2/3$. It  remains to prove that $F$ is bi-Lipschitz
and compute its distorsion.

$$\begin{aligned}||F(x)-F(y)||&\le ||F(x)-x||+||F(y)-y||+||x-y||\le 2\widetilde b+||x-y||\\&\le\bigg(\frac{2\widetilde b}{\widetilde a}+1\bigg)||x-y||\;\;\forall x,y\in N,\end{aligned}$$
and
$$\begin{aligned}||F^{-1}(x)-F^{-1}(y)||&\le ||F^{-1}(x)-x||+||F^{-1}(y)-y||+||x-y||\\&\le 2\widetilde b+||x-y||\le\bigg(\frac{4\widetilde b}{\widetilde a}+1\bigg)||x-y||\;\;\forall x,y\in S.\end{aligned}$$
We are done since $\frac{\widetilde{b}}{\widetilde{a}}=\frac{b}{a}$.
\end{proof}
\end{lemma}

From now on we are going to deal with a fixed $(a,b)$-spiderweb $S$ throughout all this section, with $a,b\in\R^+$ arbitrary. By assumption there is a base of an $(a,b)$-spiderweb $N\subset S$. We will denote $S_n=\big(nS_X\big)\cap N$ and $B_n=\big(nB_X\big)\cap S$ for every $n\in\N\cup\{0\}$. The set $B_n$ is a discrete ball of radius $n$ onto which we want to retract $S$.
Let $\rho_n:nS_X\to \frac{n}{n-1}S_{n-1}$ be the nearest point map when $n\ge2$ and  let $\rho_1:S_X\to \{0\}$ be constant.

For each $n\in\N$ we define the 'local' retraction $\psi_n:B_n\rightarrow B_{n-1}$ as
$$\psi_n(x)=\begin{cases} \frac{n-1}{n}\rho_{n}\Big(n\frac{x}{||x||}\Big)\;\;\;&\text{if }x\in B_n\setminus B_{n-1},\\x&\text{if }x\in B_{n-1}. \end{cases}$$

 Then, we define $n(x)=\min\{n\in\N\;:\;||x||\le n\}$ for every $x\in S$. Finally, the global retractions are defined for every $n\in\N\cup\{0\}$ as $\Psi_n:S\to B_n$ given by
$$\Psi_n(x)=\begin{cases}\psi_{n+1}\circ\cdots\circ\psi_{n(x)}(x) \;\;\;&\text{if }x\in S\setminus B_{n},\\x&\text{if }x\in B_{n}. \end{cases}$$

\begin{figure}\label{spiderfig}
\centering

\tikzset{every picture/.style={line width=0.75pt}} 

\begin{tikzpicture}[x=0.75pt,y=0.75pt,yscale=-1,xscale=1]

\draw  [draw opacity=0][dash pattern={on 0.84pt off 2.51pt}] (281.39,133.99) .. controls (291.83,125.24) and (305.3,119.97) .. (320,119.97) .. controls (333.94,119.97) and (346.78,124.71) .. (356.98,132.67) -- (320,180) -- cycle ; \draw  [dash pattern={on 0.84pt off 2.51pt}] (281.39,133.99) .. controls (291.83,125.24) and (305.3,119.97) .. (320,119.97) .. controls (333.94,119.97) and (346.78,124.71) .. (356.98,132.67) ;
\draw  [draw opacity=0][dash pattern={on 0.84pt off 2.51pt}] (261.94,110.6) .. controls (277.58,97.51) and (297.74,89.63) .. (319.75,89.63) .. controls (340.63,89.63) and (359.85,96.72) .. (375.12,108.63) -- (319.75,179.5) -- cycle ; \draw  [dash pattern={on 0.84pt off 2.51pt}] (261.94,110.6) .. controls (277.58,97.51) and (297.74,89.63) .. (319.75,89.63) .. controls (340.63,89.63) and (359.85,96.72) .. (375.12,108.63) ;
\draw  [draw opacity=0][dash pattern={on 0.84pt off 2.51pt}] (242.84,87.79) .. controls (263.71,70.29) and (290.62,59.75) .. (320,59.75) .. controls (347.87,59.75) and (373.53,69.24) .. (393.91,85.15) -- (320,179.75) -- cycle ; \draw  [dash pattern={on 0.84pt off 2.51pt}] (242.84,87.79) .. controls (263.71,70.29) and (290.62,59.75) .. (320,59.75) .. controls (347.87,59.75) and (373.53,69.24) .. (393.91,85.15) ;
\draw  [fill={rgb, 255:red, 0; green, 0; blue, 0 }  ,fill opacity=1 ] (293.7,93.15) .. controls (293.7,93.92) and (294.33,94.55) .. (295.1,94.55) .. controls (295.87,94.55) and (296.5,93.92) .. (296.5,93.15) .. controls (296.5,92.38) and (295.87,91.75) .. (295.1,91.75) .. controls (294.33,91.75) and (293.7,92.38) .. (293.7,93.15) -- cycle ;
\draw  [fill={rgb, 255:red, 0; green, 0; blue, 0 }  ,fill opacity=1 ] (301.95,122.15) .. controls (301.95,122.92) and (302.58,123.55) .. (303.35,123.55) .. controls (304.12,123.55) and (304.75,122.92) .. (304.75,122.15) .. controls (304.75,121.38) and (304.12,120.75) .. (303.35,120.75) .. controls (302.58,120.75) and (301.95,121.38) .. (301.95,122.15) -- cycle ;
\draw  [fill={rgb, 255:red, 0; green, 0; blue, 0 }  ,fill opacity=1 ] (335.2,122.15) .. controls (335.2,122.92) and (335.83,123.55) .. (336.6,123.55) .. controls (337.37,123.55) and (338,122.92) .. (338,122.15) .. controls (338,121.38) and (337.37,120.75) .. (336.6,120.75) .. controls (335.83,120.75) and (335.2,121.38) .. (335.2,122.15) -- cycle ;
\draw  [fill={rgb, 255:red, 0; green, 0; blue, 0 }  ,fill opacity=1 ] (343.2,93.15) .. controls (343.2,93.92) and (343.83,94.55) .. (344.6,94.55) .. controls (345.37,94.55) and (346,93.92) .. (346,93.15) .. controls (346,92.38) and (345.37,91.75) .. (344.6,91.75) .. controls (343.83,91.75) and (343.2,92.38) .. (343.2,93.15) -- cycle ;
\draw  [fill={rgb, 255:red, 0; green, 0; blue, 0 }  ,fill opacity=1 ] (351.7,64.15) .. controls (351.7,64.92) and (352.33,65.55) .. (353.1,65.55) .. controls (353.87,65.55) and (354.5,64.92) .. (354.5,64.15) .. controls (354.5,63.38) and (353.87,62.75) .. (353.1,62.75) .. controls (352.33,62.75) and (351.7,63.38) .. (351.7,64.15) -- cycle ;
\draw  [fill={rgb, 255:red, 0; green, 0; blue, 0 }  ,fill opacity=1 ] (285.45,64.15) .. controls (285.45,64.92) and (286.08,65.55) .. (286.85,65.55) .. controls (287.62,65.55) and (288.25,64.92) .. (288.25,64.15) .. controls (288.25,63.38) and (287.62,62.75) .. (286.85,62.75) .. controls (286.08,62.75) and (285.45,63.38) .. (285.45,64.15) -- cycle ;
\draw  [fill={rgb, 255:red, 0; green, 0; blue, 0 }  ,fill opacity=1 ] (321.45,59.65) .. controls (321.45,60.42) and (322.08,61.05) .. (322.85,61.05) .. controls (323.62,61.05) and (324.25,60.42) .. (324.25,59.65) .. controls (324.25,58.88) and (323.62,58.25) .. (322.85,58.25) .. controls (322.08,58.25) and (321.45,58.88) .. (321.45,59.65) -- cycle ;
\draw  [fill={rgb, 255:red, 0; green, 0; blue, 0 }  ,fill opacity=1 ] (379.95,76.65) .. controls (379.95,77.42) and (380.58,78.05) .. (381.35,78.05) .. controls (382.12,78.05) and (382.75,77.42) .. (382.75,76.65) .. controls (382.75,75.88) and (382.12,75.25) .. (381.35,75.25) .. controls (380.58,75.25) and (379.95,75.88) .. (379.95,76.65) -- cycle ;
\draw  [fill={rgb, 255:red, 0; green, 0; blue, 0 }  ,fill opacity=1 ] (254.95,77.9) .. controls (254.95,78.67) and (255.58,79.3) .. (256.35,79.3) .. controls (257.12,79.3) and (257.75,78.67) .. (257.75,77.9) .. controls (257.75,77.13) and (257.12,76.5) .. (256.35,76.5) .. controls (255.58,76.5) and (254.95,77.13) .. (254.95,77.9) -- cycle ;
\draw    (257.9,74.5) .. controls (261.77,66.77) and (270.4,61.93) .. (280.21,62.82) ;
\draw [shift={(283.15,63.25)}, rotate = 191.56] [fill={rgb, 255:red, 0; green, 0; blue, 0 }  ][line width=0.08]  [draw opacity=0] (7.14,-3.43) -- (0,0) -- (7.14,3.43) -- (4.74,0) -- cycle    ;
\draw    (326.15,57.5) .. controls (332.11,52.83) and (342.12,54.49) .. (348.47,59.32) ;
\draw [shift={(350.65,61.25)}, rotate = 226.17] [fill={rgb, 255:red, 0; green, 0; blue, 0 }  ][line width=0.08]  [draw opacity=0] (7.14,-3.43) -- (0,0) -- (7.14,3.43) -- (4.74,0) -- cycle    ;
\draw    (380.1,72.8) .. controls (377.42,65.52) and (369.89,61.76) .. (360.26,62.65) ;
\draw [shift={(357.35,63.05)}, rotate = 349.7] [fill={rgb, 255:red, 0; green, 0; blue, 0 }  ][line width=0.08]  [draw opacity=0] (7.14,-3.43) -- (0,0) -- (7.14,3.43) -- (4.74,0) -- cycle    ;
\draw    (296.53,97.4) -- (301.71,115.52) ;
\draw [shift={(302.53,118.4)}, rotate = 254.05] [fill={rgb, 255:red, 0; green, 0; blue, 0 }  ][line width=0.08]  [draw opacity=0] (7.14,-3.43) -- (0,0) -- (7.14,3.43) -- (4.74,0) -- cycle    ;
\draw    (288.2,68.07) -- (293.38,86.18) ;
\draw [shift={(294.2,89.07)}, rotate = 254.05] [fill={rgb, 255:red, 0; green, 0; blue, 0 }  ][line width=0.08]  [draw opacity=0] (7.14,-3.43) -- (0,0) -- (7.14,3.43) -- (4.74,0) -- cycle    ;
\draw    (351.53,68.87) -- (347,85.31) ;
\draw [shift={(346.2,88.2)}, rotate = 285.42] [fill={rgb, 255:red, 0; green, 0; blue, 0 }  ][line width=0.08]  [draw opacity=0] (7.14,-3.43) -- (0,0) -- (7.14,3.43) -- (4.74,0) -- cycle    ;
\draw    (343.53,97.87) -- (339,114.31) ;
\draw [shift={(338.2,117.2)}, rotate = 285.42] [fill={rgb, 255:red, 0; green, 0; blue, 0 }  ][line width=0.08]  [draw opacity=0] (7.14,-3.43) -- (0,0) -- (7.14,3.43) -- (4.74,0) -- cycle    ;
\draw  [draw opacity=0][dash pattern={on 0.84pt off 2.51pt}] (300.08,155.97) .. controls (305.42,150.78) and (312.49,147.62) .. (320.25,147.62) .. controls (327.63,147.62) and (334.38,150.48) .. (339.61,155.21) -- (320.25,180) -- cycle ; \draw  [dash pattern={on 0.84pt off 2.51pt}] (300.08,155.97) .. controls (305.42,150.78) and (312.49,147.62) .. (320.25,147.62) .. controls (327.63,147.62) and (334.38,150.48) .. (339.61,155.21) ;
\draw  [fill={rgb, 255:red, 0; green, 0; blue, 0 }  ,fill opacity=1 ] (327.83,149.07) .. controls (327.83,149.84) and (328.46,150.47) .. (329.23,150.47) .. controls (330.01,150.47) and (330.63,149.84) .. (330.63,149.07) .. controls (330.63,148.29) and (330.01,147.67) .. (329.23,147.67) .. controls (328.46,147.67) and (327.83,148.29) .. (327.83,149.07) -- cycle ;
\draw    (335.7,126.2) -- (331.16,142.64) ;
\draw [shift={(330.37,145.53)}, rotate = 285.42] [fill={rgb, 255:red, 0; green, 0; blue, 0 }  ][line width=0.08]  [draw opacity=0] (7.14,-3.43) -- (0,0) -- (7.14,3.43) -- (4.74,0) -- cycle    ;
\draw    (306.15,119.3) .. controls (312.75,112.25) and (323.69,111.94) .. (331.83,117.67) ;
\draw [shift={(334.15,119.55)}, rotate = 222.65] [fill={rgb, 255:red, 0; green, 0; blue, 0 }  ][line width=0.08]  [draw opacity=0] (7.14,-3.43) -- (0,0) -- (7.14,3.43) -- (4.74,0) -- cycle    ;

\draw (360.75,119.65) node [anchor=north west][inner sep=0.75pt]    {$C_{1}^{y}$};
\draw (397.5,69.15) node [anchor=north west][inner sep=0.75pt]    {$C_{2}^{y}$};
\draw (341.73,148.07) node [anchor=north west][inner sep=0.75pt]  [font=\normalsize]  {$C_{0}^{y} =\{y\}$};
\draw (316,147.4) node [anchor=north west][inner sep=0.75pt]    {$y$};
\draw (232,115.4) node [anchor=north west][inner sep=0.75pt]    {$\Psi _{1}$};

\end{tikzpicture}
\caption{}
\end{figure}

See figure \ref{spiderfig} for an example on how $\Psi$ works, where the picture refers to the different concentric spheres in $X$ and the points of the base of the spiderweb only.

It follows from the definition that $\Psi_n\circ\Psi_m=\Psi_{\min(n,m)}$ for every $n,m\in\N\cup\{0\}$.

\begin{theorem}\label{radial}
If $n\in\N$ then
\begin{equation}\label{eqradial}||\Psi_n(x)-R_n(x)||\le 6b\;\;\;\;\;\forall x\in S,\end{equation}
where $R_n:S\to nB_X$ is the radial projection.
\begin{proof}
Given a point $y\in S_{2^m}$ for some $m\in\N\cup\{0\}$, we want to locate the points from $ \bigcup\limits_{n=0}^\infty S_{2^n}$ that are mapped into $\{y\}$ by $\Psi_{2^m}$. To this end, we first define $C_0^y=\{y\}$ and proceed inductively for $k\ge1$, letting
$$C_k^y=\big\{x\in S_{2^{m+k}}\;:\;\rho_{2^{m+k}}(x)\in 2 C_{k-1}^y\big\}.$$
One may see in figure \ref{spiderfig} a particular example where the sets $C_k^y$ are related to the behaviour of $\psi_{2^0}$ for  some $y\in S_{2^0}$.

We claim that for every $k\ge1$,
\begin{equation}\label{mainpre}S_{2^{m+k}}\cap \Psi_{2^m}^{-1}(y)= C_{k}^y.\end{equation}
We prove \eqref{mainpre} by induction in $k$. For $k=1$, it follows from the fact that for every $x\in S_{2^{m+1}}$,
$$\Psi_{2^m}(x)=\psi_{2^m+1}\circ\cdots\circ\psi_{2^{m+1}}(x)=\bigg(\prod\limits_{j=2^m+1}^{2^{m+1}}\frac{j-1}{j}\bigg)\rho_{2^{m+1}}(x)=\frac{1}{2}\rho_{2^{m+1}}(x).$$ 
For $k\ge2$  by inductive assumption we have that $C_{k-1}^y=S_{2^{m+k-1}}\cap \Psi_{2^m}^{-1}(y)$. Since $\Psi_{2^m}=\Psi_{2^m}\circ\Psi_{2^{m+k-1}}$, we have for $x\in S_{2^{m+k}}$,
$$x\in \Psi_{2^m}^{-1}(y)\;\Leftrightarrow\;\Psi_{2^{m+k-1}}(x)\in\Psi_{2^m}^{-1}(y)\;\Leftrightarrow\;\Psi_{2^{m+k-1}}(x)\in C_{k-1}^y.$$
So it is enough to see that $\Psi_{2^{m+k-1}}(x)\in C_{k-1}^y$ if and only if $x\in C_k^y$. Indeed, 
$$\begin{aligned}\Psi_{2^{m+k-1}}(x)=&\psi_{2^{m+k-1}+1}\circ\cdots\circ\psi_{2^{m+k}}(x)\\=&\bigg(\prod\limits_{j=2^{m+k-1}+1}^{2^{m+k}}\frac{j-1}{j}\bigg)\rho_{2^{m+k}}(x)=\frac{1}{2}\rho_{2^{m+k}}(x),\end{aligned}$$
which finally proves \eqref{mainpre}. Denoting $C^y=\bigcup\limits_{k=0}^\infty C_k^y$, the following is satisfied.
\begin{equation}\label{eqpreimage}C^y= \bigg(\bigcup\limits_{n=1}^\infty S_{2^n}\bigg)\cap \Psi_{2^m}^{-1}(y).\end{equation}
Our goal now is to study the properties of $C^y$  in order to control the preimages of a point $y\in S_{2^m}$. In particular, if $F^y=\{x\in S_{||y||}\;:\;||y-x||\le 2b\}$ and $D^y=[1,\infty)F^y$, we are going to show that
\begin{equation}\label{claimmain}C^y\subset D^y\;\;\;\;\;\;\forall y\in \bigcup\limits_{n=0}^\infty S_{2^n}.\end{equation}
To prove \eqref{claimmain} we first prove the following equation for $k\in\N\cup\{0\}$.
\begin{equation}\label{contain}C_{k+1}^y\subset 2C_k^y+2bB_X\;\;\;\forall y\in \bigcup\limits_{n=0}^\infty S_{2^n}.\end{equation}
We proceed by taking the complements in \eqref{contain} for a fixed $y\in \bigcup\limits_{n=0}^\infty S_{2^n}$. If $x\in S_{||y||2^{k+1}}\setminus\big( 2C_k^y+2bB_X \big)$ then $d(x,2C_k^y)>2b$.  As $2S_{||y||2^{k}}$ is $2b$-dense in $||y||2^{k+1}S_X$, there exists a $z\in 2S_{||y||2^{k}}\setminus 2C^k_y$ such that $||x-z||\le 2b$ and therefore $\rho_{||y||2^{k+1}}(x)\notin 2C_k^y$. This means that $x\in S_{||y||2^{k+1}}\setminus C_{k+1}^y$ and \eqref{contain} is proved.

We return to the proof of \eqref{claimmain}. Our goal is to prove that $C^y\subset D^y$ for $y\in\bigcup\limits_{n\in\N}S_{2^n}$. To this end we are going to see that $\big(C_k^y+ 2bB_X\big)\cap S_{||y||2^k}\subset 2^kF^y$ for every $k\in\N\cup\{0\}$. We proceed by induction in $k$.  For $k=0$  obviously $\big(C_0^y+2bB_X\big)\cap S_{||y||}=F^y$. Let us consider the general case when $k+1\in\N$, assuming by induction that $\big(C_k^y+ 2bB_X\big)\cap S_{||y||2^k}\subset 2^kF^y$. Using \eqref{contain} we have
$$\begin{aligned}\big(C_{k+1}^y+2bB_X\big)\cap S_{||y||2^{k+1}} &\subset \big(2C_k^y+4bB_X\big)\cap S_{||y||2^{k+1}}\\& \subset 2\Big(\big( C_k^y+2bB_X\big)\cap S_{||y||2^{k}}\Big)\subset 2^{k+1}F^y,\end{aligned}$$
which ends the proof of \eqref{claimmain}.

We are now ready to prove the statement of the Lemma in two steps. First, we prove \eqref{eqradial} for points in $N$. Let us consider $n\in\N$ and $x\in N\setminus B_n$ (clearly the statement of the theorem is trivial for $x\in B_n$). Then there are $k_1,k_2\in\N\cup\{0\}$ with $k_1\le k_2$ such that
$$n\in\{2^{k_1},\dots,2^{k_1+1}-1\}\;\;\;,\;\;\;||x||\in\{2^{k_2},\dots,2^{k_2+1}-1\}.$$
Let  $y=\Psi_{2^{k_1}}(x)$ and $z= R_{2^{k_2}}(x)$.  From \eqref{eqpreimage} and  $\Psi_{2^{k_1}}(z)=\Psi_{2^{k_1}}(x)=y$ we see that $z\in C^y$.
So thanks to \eqref{contain}, we may conclude that $z\in D^y$. Then, $R_{2^{k_1}}(x)=R_{2^{k_1}}(z)=R_{||y||}(z)\in R_{||y||}(D^y)=F^y$ which means that
$$||\Psi_{2^{k_1}}(x)-R_{2^{k_1}}(x)||=||y-R_{2^{k_1}}(x)||\le 2b.$$
Consequently,
$$||\Psi_n(x)-R_n(x)||=\frac{n}{2^{k_1}}||\Psi_{2^{k_1}}(x)-R_{2^{k_1}}(x)||< 4b,$$
so we conclude that
\begin{equation}\label{basespider}||\Psi_n(x)-R_n(x)||\le 4b\;\;\;\;\;\;\forall x\in N.\end{equation}
Next, we consider the general case when $x\in S\setminus B_n$. Then $n(x)\ge n$ and
$$\begin{aligned}\bigg|\bigg|R_n\bigg( \rho_{n(x)}\bigg( \frac{n(x)x}{||x||} \bigg) \bigg)-R_n(x)\bigg|\bigg|&=\bigg|\bigg|\frac{n}{n(x)}\rho_{n(x)}\bigg(\frac{n(x)x}{||x||} \bigg)-\frac{n}{n(x)}\bigg(\frac{n(x)x}{||x||}\bigg)\bigg|\bigg|\\&=\frac{n}{n(x)}\bigg|\bigg| \rho_{n(x)}\bigg(\frac{n(x)x}{||x||}\bigg)-\frac{n(x)x}{||x||} \bigg|\bigg|\le2b.\end{aligned}$$
By definition we have that $\Psi_n(x)=\Psi_n\big(\rho_{n(x)}\big( \frac{n(x)x}{||x||} \big)\big)$ where $\frac{n(x)x}{||x||}\in N$. So taking into account \eqref{basespider} and the previous inequality we have that
$$\begin{aligned}||\Psi_n(x)-R_n(x)||\le&\bigg|\bigg|\Psi_n\bigg(\rho_{n(x)}\bigg( \frac{n(x)x}{||x||} \bigg)\bigg)-R_n\bigg(\rho_{n(x)}\bigg( \frac{n(x)x}{||x||} \bigg)\bigg)\bigg|\bigg|\\&+\bigg|\bigg| R_n\bigg(\rho_{n(x)}\bigg( \frac{n(x)x}{||x||} \bigg)\bigg)-R_n(x) \bigg|\bigg|\le 4b+2b.\end{aligned}$$
\end{proof}
\end{theorem}

\begin{corollary}\label{theospider}
Every $(a,b)$-spiderweb $S$ of a finite dimensional Banach space $X$ has a $K$-retractional basis where
$$K=\frac{12b+2}{a}+2.$$
\begin{proof}
Denote $B_{-1}=\emptyset$ and $k_n=\#(B_n\setminus B_{n-1})$ for every $n\in\N\cup\{0\}$. Recall that
$$S=\{x_{(n,i)}\}_{\substack{n\in\N\cup\{0\}\\i=1,\dots,k_n}},$$
where $\{x_{(n,1)},\dots,x_{(n,k_n)}\}=B_{n}\setminus B_{n-1}$ is an arbitrary order of $B_n\setminus B_{n-1}$ for every $n\in\N\cup\{0\}$.
Order S lexicographically according to the previous notation and define the retractions of the retractional basis $\varphi_{(n,i)}:S\to S_{(n,i)}:=\{x_{(m,j)}\}_{(m,j)\le(n,i)}$ as
$$\varphi_{(n,i)}(x)=\begin{cases}\Psi_n(x)\;\;&\text{if }\Psi_n(x)\in S_{(n,i)},\\\Psi_{n-1}(x)\;\;&\text{if }\Psi_n(x)\notin S_{(n,i)}.\end{cases}$$
It is easy to see that $\varphi_{(n,i)}$ is a well defined retraction. Also, a straightforward computation gives $\varphi_{(n,i)}\circ\varphi_{(m,j)}=\varphi_{\min((n,i),(m,j))}$ so that it only remains to prove that $\varphi_{(n,i)}$ is Lipschitz. Indeed, by Lemma \ref{radial}, if $x\in S$ then
$$||\varphi_{(n,i)}(x)-R_{n}(x)||\le\max\{||\Psi_n(x)-R_n(x)||,||\Psi_{n-1}(x)-R_n(x)||\}\le 6b+1.$$
Then for every $x,y\in S$,
$$||\varphi_{(n,i)}(x)-\varphi_{(n,i)}(y)||\le 2(6b+1)+||R_n(x)-R_n(y)||\le \bigg(\frac{12b+2}{a}+2\bigg)||x-y||.$$
\end{proof}
\end{corollary}

\begin{theorem}\label{theofinite}
Every $(a,b)$-net of a finite dimensional Banach space $X$ has a $K$-retractional basis where 
$$K=\bigg( \frac{2b}{a}+1 \bigg)\bigg( \frac{4b}{a}+1 \bigg)\bigg( \frac{156b}{a}+2 \bigg).$$
\begin{proof}
Let $N$ be an arbitrary $(a,b)$-net of $X$. By Lemma \ref{spider}  $N$ is Lipschitz equivalent  to a $\big(\frac{a}{6b},2\big)$-spiderweb $S$ of $X$, with distorsion $D=\big(\frac{2b}{a}+1\big)\big(\frac{4b}{a}+1\big)$. By corollary \ref{theospider}, S has a $\widetilde{K}$-retractional basis where
$$\widetilde{K}=\frac{12\cdot 2+2}{\frac{a}{6b}}+2=\frac{156b}{a}+2.$$
Finally, thanks to Proposition \ref{lipequiv} we know that $N$ has a $K$-retractional basis where
$$K=D\widetilde K =\bigg( \frac{2b}{a}+1 \bigg)\bigg( \frac{4b}{a}+1 \bigg)\bigg( \frac{156b}{a}+2 \bigg).$$
\end{proof}
\end{theorem}

\begin{remark}
It is worth  mentioning that $b/a$ may serve as a measure of  homogeneity of an $(a,b)$-net. We have not strived to obtain an optimal
value of $K$ in the previous theorem.
\end{remark}

\section{FDD grids}\label{FDD}


In the present section we will study the retractional basis when $N$ is a certain grid in a Banach space $X$ admitting an FDD. 
In general, such grids are not necessarily nets in $X$ (i.e. they may not be sufficiently dense in $X$), unless $X$ has some additional
$c_0$-like structure. Our retractional technique is a spin-off of our methods in \cite{HM21}. In the mentioned paper we have
constructed a diamond shaped generating convex compact set in $X$, with fast decreasing sequence of  "heights", which served
as a target for a Lipschitz retraction from $X$. In the present paper we proceed on the contrary by blowing up the generating diamond, creating
an inverse-limit type of situation based on diamond-like (but discrete, in fact finite) subsets of the grid of a growing size. 
The diamond shapes are precisely tuned to match the grid in $X$, in the sense that every element of the grid
is in some sense a "boundary" point of a certain diamond and then it becomes an "interior" point for all larger diamonds. The main mechanism of the construction
of our retractions consists roughly speaking of reducing the last non-zero coordinate (whose norm is  an integer) of an element of the grid to a coordinate of  norm smaller by one. In this way,
the element is moved from its boundary position, with respect to a specific diamond, to an interior position thereof.  
The crucial ingredient is also the independence
of the retractional constant from the previous finite dimensional case on the dimension.

\begin{definition}
Let $X$ be a Banach space with an FDD $(X_n)$ and $a,b\in\R^+$. A subset $S\subset X$ is an $(a,b)$-spiderweb grid with respect to the FDD if
$$S=\bigcup\limits_{n\in\N}\sum\limits_{m=1}^nS^m,$$
where $S^m$ is a base of an $(a,b)$-spiderweb of $X_m$, for every $m\in\N$.
\end{definition}

From now on $X$ will be a Banach space with an FDD $(X_n)$ and
$S=\bigcup\limits_{n\in\N}\sum\limits_{m=1}^nS^m$
an $(a,b)$-spiderweb grid with respect to this FDD. We denote the $n^{th}$ coordinate of $x\in X$ by $x_n=\big(P_n-P_{n-1}\big)(x)$, where $(P_n)$ is the sequence of the natural projections of the FDD. We proceed by  defining the diamond-shaped sets needed for our construction.

For every $s\in \N\cup\{0\}$ we let
$$D^s=\cco\bigg( \bigcup\limits_{k\in\N}r^s_kB_ {X_k}\bigg),$$
where the sequence $(r_k^s)_{k\in\N}$ satisfies
$$\begin{cases}
r_1^s=s,\\
r^s_{k+1}=\frac{1}{k2^{k+2}}r^s_k.
\end{cases}$$
Notice that $D^0=\{0\}$ and every $D^s$ is compact since $r_n^s\xrightarrow{n\to\infty}0$.
Throughout  we will agree that $\sum\limits_{i=n}^{m}a_i=0$ for every sequence $(a_i)$ and $n,m\in\Z$ such that $n<m$. From the definition of the sequences $(r_n^s)_{n\in\N}$ it follows easily that

\begin{equation}\label{proprn}\frac{r^s_n}{r^s_m}=\frac{r^t_n}{r^t_m}\in\N\;,\;\;\;\;\;\frac{r^s_n}{r^t_n}=\frac{s}{t}\;,\;\;\;\;\;\;\;\;\;\;\;\forall n,m,s,t\in \N\;,\;\;n\le m.\end{equation}

We begin the construction of the retractional basis for $S$ by defining an increasing sequence of finite subsets $(N_s)_{s\in\N\cup\{0\}}$ of $S$ as
$$N_s=D^s\cap S\;\;\;\;\forall s\in \N \cup\{0\}.$$
The set $N_s$ is a discrete diamond whose highest length is $s\in\N\cup\{0\}$. Obviously, $N_0=\{0\}$ and for aesthetic purposes we set $N_{-1}=\emptyset$. For each $n\in\N$ we consider the sequence of retractions $(\Psi^n_s)_{s\in\N\cup\{0\}}$ where $\Psi^n_s:S^n\rightarrow sB_{X_n}\cap S^n$ is the retraction defined in Section \ref{finite}, but now the radius of the discrete ball onto which we retract is given by $s\in\N\cup\{0\}$ (in section \ref{finite} we used $n$). Recall that
\begin{equation}\label{similrad}||\Psi^n_s(x)-R_s(x)||\le 6b\;\;\;\;\forall x\in S^n,\end{equation}
where $R_s:X\rightarrow sB_{X}$ is the radial projection. To simplify the notation we will abbreviate $\Psi^n_s$  as $\Psi_s$ since the superindex $n$ only indicates the changing domain. For every $x\in \bigcup\limits_{n\in\N}P_n(X)\setminus\{0\}$, we are going to denote $n(x)=\max\{n\in\N\;:\;x_n\neq0\}$. Next, we define the 'local' retractions $\varphi_s:N_s\to N_{s-1}$ as
$$\varphi_s(x)=\begin{cases} x\;\;\;&\text{if }x\in N_{s-1},\\P_{n(x)-1}(x)+\Psi_{||x_{n(x)}||-1}(x_{n(x)})\;\;\;&\text{if }x\in N_s\setminus N_{s-1}. \end{cases}$$
For each $x\in S$ let us define $s(x)$ as the unique $s\in \N\cup\{0\}$ such that $x\in N_{s}\setminus N_{s-1}$.
Finally, we are ready to define the retractions $\phi_s:S\to N_{s}$ for every $s\in\N\cup\{0\}$ as
$$\phi_s(x)=\begin{cases}\varphi_{s+1}\circ\cdots\circ\varphi_{s(x)}(x)\;\;\;&\text{if }x\in S\setminus N_s,\\x&\text{if }x\in N_s.\end{cases}$$

\begin{figure}\label{diamond1}
\centering

\tikzset{every picture/.style={line width=0.75pt}} 

\begin{tikzpicture}[x=0.75pt,y=0.75pt,yscale=-1,xscale=1]

\draw  [fill={rgb, 255:red, 0; green, 0; blue, 0 }  ,fill opacity=1 ] (105,152.5) .. controls (105,151.12) and (106.12,150) .. (107.5,150) .. controls (108.88,150) and (110,151.12) .. (110,152.5) .. controls (110,153.88) and (108.88,155) .. (107.5,155) .. controls (106.12,155) and (105,153.88) .. (105,152.5) -- cycle ;
\draw  [fill={rgb, 255:red, 0; green, 0; blue, 0 }  ,fill opacity=1 ] (105,182.5) .. controls (105,181.12) and (106.12,180) .. (107.5,180) .. controls (108.88,180) and (110,181.12) .. (110,182.5) .. controls (110,183.88) and (108.88,185) .. (107.5,185) .. controls (106.12,185) and (105,183.88) .. (105,182.5) -- cycle ;
\draw  [fill={rgb, 255:red, 0; green, 0; blue, 0 }  ,fill opacity=1 ] (135,152.5) .. controls (135,151.12) and (136.12,150) .. (137.5,150) .. controls (138.88,150) and (140,151.12) .. (140,152.5) .. controls (140,153.88) and (138.88,155) .. (137.5,155) .. controls (136.12,155) and (135,153.88) .. (135,152.5) -- cycle ;
\draw  [fill={rgb, 255:red, 0; green, 0; blue, 0 }  ,fill opacity=1 ] (135,182.5) .. controls (135,181.12) and (136.12,180) .. (137.5,180) .. controls (138.88,180) and (140,181.12) .. (140,182.5) .. controls (140,183.88) and (138.88,185) .. (137.5,185) .. controls (136.12,185) and (135,183.88) .. (135,182.5) -- cycle ;
\draw  [fill={rgb, 255:red, 0; green, 0; blue, 0 }  ,fill opacity=1 ] (165,152.5) .. controls (165,151.12) and (166.12,150) .. (167.5,150) .. controls (168.88,150) and (170,151.12) .. (170,152.5) .. controls (170,153.88) and (168.88,155) .. (167.5,155) .. controls (166.12,155) and (165,153.88) .. (165,152.5) -- cycle ;
\draw  [fill={rgb, 255:red, 0; green, 0; blue, 0 }  ,fill opacity=1 ] (165,182.5) .. controls (165,181.12) and (166.12,180) .. (167.5,180) .. controls (168.88,180) and (170,181.12) .. (170,182.5) .. controls (170,183.88) and (168.88,185) .. (167.5,185) .. controls (166.12,185) and (165,183.88) .. (165,182.5) -- cycle ;
\draw  [fill={rgb, 255:red, 0; green, 0; blue, 0 }  ,fill opacity=1 ] (195,182.5) .. controls (195,181.12) and (196.12,180) .. (197.5,180) .. controls (198.88,180) and (200,181.12) .. (200,182.5) .. controls (200,183.88) and (198.88,185) .. (197.5,185) .. controls (196.12,185) and (195,183.88) .. (195,182.5) -- cycle ;
\draw  [fill={rgb, 255:red, 0; green, 0; blue, 0 }  ,fill opacity=1 ] (225,182.5) .. controls (225,181.12) and (226.12,180) .. (227.5,180) .. controls (228.88,180) and (230,181.12) .. (230,182.5) .. controls (230,183.88) and (228.88,185) .. (227.5,185) .. controls (226.12,185) and (225,183.88) .. (225,182.5) -- cycle ;
\draw  [fill={rgb, 255:red, 0; green, 0; blue, 0 }  ,fill opacity=1 ] (105,122.5) .. controls (105,121.12) and (106.12,120) .. (107.5,120) .. controls (108.88,120) and (110,121.12) .. (110,122.5) .. controls (110,123.88) and (108.88,125) .. (107.5,125) .. controls (106.12,125) and (105,123.88) .. (105,122.5) -- cycle ;
\draw   (135,122.5) .. controls (135,121.12) and (136.12,120) .. (137.5,120) .. controls (138.88,120) and (140,121.12) .. (140,122.5) .. controls (140,123.88) and (138.88,125) .. (137.5,125) .. controls (136.12,125) and (135,123.88) .. (135,122.5) -- cycle ;
\draw   (165,122.5) .. controls (165,121.12) and (166.12,120) .. (167.5,120) .. controls (168.88,120) and (170,121.12) .. (170,122.5) .. controls (170,123.88) and (168.88,125) .. (167.5,125) .. controls (166.12,125) and (165,123.88) .. (165,122.5) -- cycle ;
\draw   (195,122.5) .. controls (195,121.12) and (196.12,120) .. (197.5,120) .. controls (198.88,120) and (200,121.12) .. (200,122.5) .. controls (200,123.88) and (198.88,125) .. (197.5,125) .. controls (196.12,125) and (195,123.88) .. (195,122.5) -- cycle ;
\draw   (195,152.5) .. controls (195,151.12) and (196.12,150) .. (197.5,150) .. controls (198.88,150) and (200,151.12) .. (200,152.5) .. controls (200,153.88) and (198.88,155) .. (197.5,155) .. controls (196.12,155) and (195,153.88) .. (195,152.5) -- cycle ;
\draw   (225,122.5) .. controls (225,121.12) and (226.12,120) .. (227.5,120) .. controls (228.88,120) and (230,121.12) .. (230,122.5) .. controls (230,123.88) and (228.88,125) .. (227.5,125) .. controls (226.12,125) and (225,123.88) .. (225,122.5) -- cycle ;
\draw   (225,152.5) .. controls (225,151.12) and (226.12,150) .. (227.5,150) .. controls (228.88,150) and (230,151.12) .. (230,152.5) .. controls (230,153.88) and (228.88,155) .. (227.5,155) .. controls (226.12,155) and (225,153.88) .. (225,152.5) -- cycle ;
\draw    (107.5,122.5) -- (227.5,182.5) ;
\draw    (107.5,122.5) -- (92.68,129.92) ;
\draw    (227.5,182.5) -- (212.68,189.92) ;
\draw    (107.5,137.7) -- (197.5,182.5) ;
\draw    (197.5,182.5) -- (182.68,189.92) ;
\draw    (107.5,137.7) -- (92.68,145.12) ;
\draw  [dash pattern={on 0.84pt off 2.51pt}]  (122.3,114.9) -- (242.3,174.9) ;
\draw  [dash pattern={on 0.84pt off 2.51pt}]  (152.3,114.9) -- (242.28,160) ;
\draw  [dash pattern={on 0.84pt off 2.51pt}]  (182.3,114.9) -- (242.29,145.25) ;
\draw  [dash pattern={on 0.84pt off 2.51pt}]  (212.3,114.9) -- (242.68,130.4) ;
\draw  [dash pattern={on 0.84pt off 2.51pt}]  (107.5,152.5) -- (92.68,159.92) ;
\draw  [dash pattern={on 0.84pt off 2.51pt}]  (167.5,182.5) -- (152.68,189.92) ;
\draw  [dash pattern={on 0.84pt off 2.51pt}]  (107.5,152.5) -- (167.5,182.5) ;
\draw  [dash pattern={on 0.84pt off 2.51pt}]  (107.5,167.5) -- (137.5,182.5) ;
\draw  [dash pattern={on 0.84pt off 2.51pt}]  (107.5,167.5) -- (92.68,174.92) ;
\draw  [dash pattern={on 0.84pt off 2.51pt}]  (137.5,182.5) -- (122.68,189.92) ;
\draw    (107.57,126.94) -- (107.51,144.45) ;
\draw [shift={(107.5,147.45)}, rotate = 270.2] [fill={rgb, 255:red, 0; green, 0; blue, 0 }  ][line width=0.08]  [draw opacity=0] (8.04,-3.86) -- (0,0) -- (8.04,3.86) -- (5.34,0) -- cycle    ;
\draw    (167.57,156.94) -- (167.51,174.45) ;
\draw [shift={(167.5,177.45)}, rotate = 270.2] [fill={rgb, 255:red, 0; green, 0; blue, 0 }  ][line width=0.08]  [draw opacity=0] (8.04,-3.86) -- (0,0) -- (8.04,3.86) -- (5.34,0) -- cycle    ;
\draw    (221.96,182.4) -- (205.36,182.4) ;
\draw [shift={(202.36,182.4)}, rotate = 360] [fill={rgb, 255:red, 0; green, 0; blue, 0 }  ][line width=0.08]  [draw opacity=0] (8.04,-3.86) -- (0,0) -- (8.04,3.86) -- (5.34,0) -- cycle    ;
\draw  [fill={rgb, 255:red, 0; green, 0; blue, 0 }  ,fill opacity=1 ] (323,152.5) .. controls (323,151.12) and (324.12,150) .. (325.5,150) .. controls (326.88,150) and (328,151.12) .. (328,152.5) .. controls (328,153.88) and (326.88,155) .. (325.5,155) .. controls (324.12,155) and (323,153.88) .. (323,152.5) -- cycle ;
\draw  [fill={rgb, 255:red, 0; green, 0; blue, 0 }  ,fill opacity=1 ] (323,182.5) .. controls (323,181.12) and (324.12,180) .. (325.5,180) .. controls (326.88,180) and (328,181.12) .. (328,182.5) .. controls (328,183.88) and (326.88,185) .. (325.5,185) .. controls (324.12,185) and (323,183.88) .. (323,182.5) -- cycle ;
\draw  [fill={rgb, 255:red, 0; green, 0; blue, 0 }  ,fill opacity=1 ] (353,152.5) .. controls (353,151.12) and (354.12,150) .. (355.5,150) .. controls (356.88,150) and (358,151.12) .. (358,152.5) .. controls (358,153.88) and (356.88,155) .. (355.5,155) .. controls (354.12,155) and (353,153.88) .. (353,152.5) -- cycle ;
\draw  [fill={rgb, 255:red, 0; green, 0; blue, 0 }  ,fill opacity=1 ] (353,182.5) .. controls (353,181.12) and (354.12,180) .. (355.5,180) .. controls (356.88,180) and (358,181.12) .. (358,182.5) .. controls (358,183.88) and (356.88,185) .. (355.5,185) .. controls (354.12,185) and (353,183.88) .. (353,182.5) -- cycle ;
\draw   (383,152.5) .. controls (383,151.12) and (384.12,150) .. (385.5,150) .. controls (386.88,150) and (388,151.12) .. (388,152.5) .. controls (388,153.88) and (386.88,155) .. (385.5,155) .. controls (384.12,155) and (383,153.88) .. (383,152.5) -- cycle ;
\draw  [fill={rgb, 255:red, 0; green, 0; blue, 0 }  ,fill opacity=1 ] (383,182.5) .. controls (383,181.12) and (384.12,180) .. (385.5,180) .. controls (386.88,180) and (388,181.12) .. (388,182.5) .. controls (388,183.88) and (386.88,185) .. (385.5,185) .. controls (384.12,185) and (383,183.88) .. (383,182.5) -- cycle ;
\draw  [fill={rgb, 255:red, 0; green, 0; blue, 0 }  ,fill opacity=1 ] (413,182.5) .. controls (413,181.12) and (414.12,180) .. (415.5,180) .. controls (416.88,180) and (418,181.12) .. (418,182.5) .. controls (418,183.88) and (416.88,185) .. (415.5,185) .. controls (414.12,185) and (413,183.88) .. (413,182.5) -- cycle ;
\draw   (443,182.5) .. controls (443,181.12) and (444.12,180) .. (445.5,180) .. controls (446.88,180) and (448,181.12) .. (448,182.5) .. controls (448,183.88) and (446.88,185) .. (445.5,185) .. controls (444.12,185) and (443,183.88) .. (443,182.5) -- cycle ;
\draw   (323,122.5) .. controls (323,121.12) and (324.12,120) .. (325.5,120) .. controls (326.88,120) and (328,121.12) .. (328,122.5) .. controls (328,123.88) and (326.88,125) .. (325.5,125) .. controls (324.12,125) and (323,123.88) .. (323,122.5) -- cycle ;
\draw   (353,122.5) .. controls (353,121.12) and (354.12,120) .. (355.5,120) .. controls (356.88,120) and (358,121.12) .. (358,122.5) .. controls (358,123.88) and (356.88,125) .. (355.5,125) .. controls (354.12,125) and (353,123.88) .. (353,122.5) -- cycle ;
\draw   (383,122.5) .. controls (383,121.12) and (384.12,120) .. (385.5,120) .. controls (386.88,120) and (388,121.12) .. (388,122.5) .. controls (388,123.88) and (386.88,125) .. (385.5,125) .. controls (384.12,125) and (383,123.88) .. (383,122.5) -- cycle ;
\draw   (413,122.5) .. controls (413,121.12) and (414.12,120) .. (415.5,120) .. controls (416.88,120) and (418,121.12) .. (418,122.5) .. controls (418,123.88) and (416.88,125) .. (415.5,125) .. controls (414.12,125) and (413,123.88) .. (413,122.5) -- cycle ;
\draw   (413,152.5) .. controls (413,151.12) and (414.12,150) .. (415.5,150) .. controls (416.88,150) and (418,151.12) .. (418,152.5) .. controls (418,153.88) and (416.88,155) .. (415.5,155) .. controls (414.12,155) and (413,153.88) .. (413,152.5) -- cycle ;
\draw   (443,122.5) .. controls (443,121.12) and (444.12,120) .. (445.5,120) .. controls (446.88,120) and (448,121.12) .. (448,122.5) .. controls (448,123.88) and (446.88,125) .. (445.5,125) .. controls (444.12,125) and (443,123.88) .. (443,122.5) -- cycle ;
\draw   (443,152.5) .. controls (443,151.12) and (444.12,150) .. (445.5,150) .. controls (446.88,150) and (448,151.12) .. (448,152.5) .. controls (448,153.88) and (446.88,155) .. (445.5,155) .. controls (444.12,155) and (443,153.88) .. (443,152.5) -- cycle ;
\draw  [dash pattern={on 0.84pt off 2.51pt}]  (325.5,122.5) -- (386.54,153.02) -- (445.5,182.5) ;
\draw  [dash pattern={on 0.84pt off 2.51pt}]  (325.5,122.5) -- (310.68,129.92) ;
\draw  [dash pattern={on 0.84pt off 2.51pt}]  (445.5,182.5) -- (430.68,189.92) ;
\draw    (325.5,137.7) -- (415.5,182.5) ;
\draw    (415.5,182.5) -- (400.68,189.92) ;
\draw    (325.5,137.7) -- (310.68,145.12) ;
\draw  [dash pattern={on 0.84pt off 2.51pt}]  (340.3,114.9) -- (460.3,174.9) ;
\draw  [dash pattern={on 0.84pt off 2.51pt}]  (370.3,114.9) -- (460.28,160) ;
\draw  [dash pattern={on 0.84pt off 2.51pt}]  (400.3,114.9) -- (460.29,145.25) ;
\draw  [dash pattern={on 0.84pt off 2.51pt}]  (430.3,114.9) -- (460.68,130.4) ;
\draw  [dash pattern={on 0.84pt off 2.51pt}]  (325.5,152.5) -- (310.68,159.92) ;
\draw  [dash pattern={on 0.84pt off 2.51pt}]  (385.5,182.5) -- (370.68,189.92) ;
\draw  [dash pattern={on 0.84pt off 2.51pt}]  (325.5,152.5) -- (385.5,182.5) ;
\draw  [dash pattern={on 0.84pt off 2.51pt}]  (325.5,167.5) -- (355.5,182.5) ;
\draw  [dash pattern={on 0.84pt off 2.51pt}]  (325.5,167.5) -- (310.68,174.92) ;
\draw  [dash pattern={on 0.84pt off 2.51pt}]  (355.5,182.5) -- (340.68,189.92) ;
\draw    (325.57,126.94) -- (325.51,144.45) ;
\draw [shift={(325.5,147.45)}, rotate = 270.2] [fill={rgb, 255:red, 0; green, 0; blue, 0 }  ][line width=0.08]  [draw opacity=0] (8.04,-3.86) -- (0,0) -- (8.04,3.86) -- (5.34,0) -- cycle    ;
\draw    (385.57,156.94) -- (385.51,174.45) ;
\draw [shift={(385.5,177.45)}, rotate = 270.2] [fill={rgb, 255:red, 0; green, 0; blue, 0 }  ][line width=0.08]  [draw opacity=0] (8.04,-3.86) -- (0,0) -- (8.04,3.86) -- (5.34,0) -- cycle    ;
\draw    (439.96,182.4) -- (423.36,182.4) ;
\draw [shift={(420.36,182.4)}, rotate = 360] [fill={rgb, 255:red, 0; green, 0; blue, 0 }  ][line width=0.08]  [draw opacity=0] (8.04,-3.86) -- (0,0) -- (8.04,3.86) -- (5.34,0) -- cycle    ;
\draw    (355.57,127.34) -- (355.51,144.85) ;
\draw [shift={(355.5,147.85)}, rotate = 270.2] [fill={rgb, 255:red, 0; green, 0; blue, 0 }  ][line width=0.08]  [draw opacity=0] (8.04,-3.86) -- (0,0) -- (8.04,3.86) -- (5.34,0) -- cycle    ;
\draw    (385.5,127.4) -- (385.44,144.91) ;
\draw [shift={(385.43,147.91)}, rotate = 270.2] [fill={rgb, 255:red, 0; green, 0; blue, 0 }  ][line width=0.08]  [draw opacity=0] (8.04,-3.86) -- (0,0) -- (8.04,3.86) -- (5.34,0) -- cycle    ;
\draw    (415.32,127.19) -- (415.27,144.71) ;
\draw [shift={(415.25,147.71)}, rotate = 270.2] [fill={rgb, 255:red, 0; green, 0; blue, 0 }  ][line width=0.08]  [draw opacity=0] (8.04,-3.86) -- (0,0) -- (8.04,3.86) -- (5.34,0) -- cycle    ;
\draw    (445.57,127.34) -- (445.51,144.85) ;
\draw [shift={(445.5,147.85)}, rotate = 270.2] [fill={rgb, 255:red, 0; green, 0; blue, 0 }  ][line width=0.08]  [draw opacity=0] (8.04,-3.86) -- (0,0) -- (8.04,3.86) -- (5.34,0) -- cycle    ;
\draw    (415.57,156.84) -- (415.51,174.35) ;
\draw [shift={(415.5,177.35)}, rotate = 270.2] [fill={rgb, 255:red, 0; green, 0; blue, 0 }  ][line width=0.08]  [draw opacity=0] (8.04,-3.86) -- (0,0) -- (8.04,3.86) -- (5.34,0) -- cycle    ;
\draw    (445.32,156.94) -- (445.26,174.45) ;
\draw [shift={(445.25,177.45)}, rotate = 270.2] [fill={rgb, 255:red, 0; green, 0; blue, 0 }  ][line width=0.08]  [draw opacity=0] (8.04,-3.86) -- (0,0) -- (8.04,3.86) -- (5.34,0) -- cycle    ;

\draw (216.67,188.93) node [anchor=north west][inner sep=0.75pt]    {$D^{s+1}$};
\draw (54.17,112.87) node [anchor=north west][inner sep=0.75pt]    {$\varphi _{s+1}$};
\draw (186.27,188.93) node [anchor=north west][inner sep=0.75pt]    {$D^{s}$};
\draw (406.37,188.93) node [anchor=north west][inner sep=0.75pt]    {$D^{s}$};
\draw (287.77,112.37) node [anchor=north west][inner sep=0.75pt]    {$\phi _{s}$};

\end{tikzpicture}

\caption{The lines of the picture refer to the different diamond shaped sets covering the space.}
\end{figure}
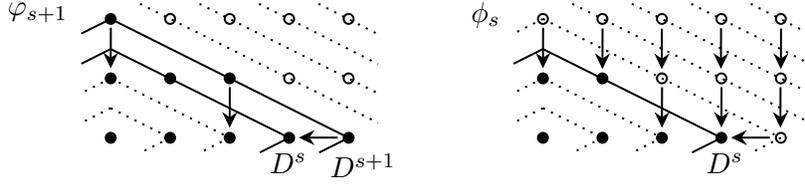
One may see in figure 3 the action of $\varphi_{s+1}$ and $\phi_s$.
Given $n\in\N$ we are going to denote $q_{n}=\frac{1}{r^1_n}=\frac{r_1^s}{r_n^s}\in\N$ and $q_0=0$.

\begin{lemma}\label{proplemma}
For every $x\in S\setminus\{0\}$ the following properties are satisfied.
\begin{enumerate}[i)]
\item $$s(x)=\sum\limits_{i=1}^{n(x)}||x_i||q_i,$$\label{sx}
\item $$\sum\limits_{i=1}^{n(x)}\frac{||x_i||}{r^{s(x)}_i}=1,$$\label{sphere}
\item $$s(\phi_{s(x)-1}(x))=s(\varphi_{s(x)}(x))=s(x)-q_{n(x)}.$$\label{iter}
\begin{proof}
Let us prove both $\ref{sx})$ and $\ref{sphere})$ together by induction in $n(x)$. If $n(x)=1$ it is trivially satisfied. We assume that it is true for every $x\in S$ such that $n(x)\le k-1$ for a fixed $k\in\N\setminus\{1\}$. Now, if $x\in S$ is such that $n(x)=k$ then $s(P_{k-1}(x))=\sum\limits_{i=1}^{k-1}||x_i||q_i$ and $\sum\limits_{i=1}^{k-1}\frac{||x_i||}{r_i^t}=1$, where $t=s(P_{k-1}(x))$. It is then enough to prove that $\ref{sphere})$ is satisfied for $s(x)=t+||x_k||q_k$. By $\eqref{proprn}$ we know that $r_k^{t+||x_k||q_k}=(t+||x_k||q_k)r_k^1$, so
$$\frac{1}{r_k^{t+||x_k||q_k}}=\frac{q_k}{t+||x_k||q_k}.$$
 If $t\neq0$ then by \eqref{proprn} and the previous equality,
$$\begin{aligned}\sum\limits_{i=1}^k\frac{||x_i||}{r_i^{t+||x_k||q_{k}}}&=\bigg(\frac{t}{t+||x_k||q_{k}}\sum\limits_{i=1}^{k-1}\frac{||x_i||}{r_i^t}\bigg)+\frac{||x_k||}{r_k^{t+||x_k||q_{k}}}\\&=\frac{t}{t+||x_k||q_{k}}+\frac{||x_k||q_{k}}{t+||x_k||q_{k}}=1.\end{aligned}$$
Otherwise, if $t=0$ then
$$\sum\limits_{i=1}^k\frac{||x_i||}{r_i^{t+||x_k||q_k}}=\frac{||x_k||}{r_k^{||x_k||+q_k}}=\frac{||x_k||q_k}{||x_k||q_k}=1,$$
finishing the proof of $\ref{sx})$ and $\ref{sphere})$.
Finally, to prove $\ref{iter})$ we let $y=\phi_{s(x)-1}(x)=\varphi_{s(x)}(x)$.
If $y=0$ then by $\ref{sx})$ we know that $s(x)=q_{n(x)}$ and trivially $y=0\in N_0=N_{s(x)-q_{n(x)}}\setminus N_{s(x)-q_{n(x)}-1}$. Otherwise, $s(x)-q_{n(x)}>0$ so similarly to the argument used in $\ref{sx})$ and $\ref{sphere})$,
$$\begin{aligned}\sum\limits_{i=1}^{n(y)}\frac{||y_i||}{r_i^{s(x)-q_{n(x)}}}&=\bigg( \sum\limits_{i=1}^{n(x)-1}                         \frac{||x_i||}{r_i^{s(x)-q_{n(x)}}} \bigg)+\frac{||x_{n(x)}||-1}{r_{n(x)}^{s(x)-q_{n(x)}}}\\&=\frac{s(x)}{s(x)-q_{n(x)}}\bigg( \sum\limits_{i=1}^{n(x)}\frac{||x_i||}{r_i^{s(x)}} \bigg)-\frac{1}{r_{n(x)}^{s(x)-q_{n(x)}}}\\&=\frac{s(x)}{s(x)-q_{n(x)}}-\frac{q_{n(x)}}{s(x)-q_{n(x)}}=1,\end{aligned}$$
and we are done.
\end{proof}
\end{enumerate}
\end{lemma}

Let us point out that $\ref{sphere})$ in Lemma \ref{proplemma} means that every point of $S\setminus\{0\}$ lies on the boundary of a diamond $D^s$ for some suitable $s\in \N$ (in fact the suitable $s$ for $x$ is $s(x)$).

It follows from the definition that
\begin{equation}\label{min}\phi_m\circ\phi_n=\phi_{min(m,n)}\;\;\;\;\;\forall m,n\in\N\cup\{0\}.\end{equation}
Every point $x\in S$ gets iteratively mapped by the local retractions when it is retracted by some $\phi_s$. Some of these local retractions act like an identity on $x$ and others  (which we will call nontrivial local retractions) do not. In order to obtain the precise information on how a point gets mapped throughout the nontrivial local retractions we define for each $x\in S$ the sequence $(s_n(x))_{n=0}^\infty$ by
$$\begin{cases} s_0(x)=s(x),\\s_{n}(x)=s(\phi_{s_{n-1}(x)-1}(x))\;\;\;\forall n\in\N,\;\end{cases}$$
where we set $\phi_{-1}:S\to \{0\}$ to be constant. 

\begin{figure}\label{IJ1}
\centering

\tikzset{every picture/.style={line width=0.75pt}} 

\begin{tikzpicture}[x=0.75pt,y=0.75pt,yscale=-1,xscale=1]

\draw    (26.53,174.29) -- (26.53,56.76) ;
\draw [shift={(26.53,54.76)}, rotate = 450] [color={rgb, 255:red, 0; green, 0; blue, 0 }  ][line width=0.75]    (10.93,-3.29) .. controls (6.95,-1.4) and (3.31,-0.3) .. (0,0) .. controls (3.31,0.3) and (6.95,1.4) .. (10.93,3.29)   ;
\draw [shift={(26.53,114.52)}, rotate = 450] [color={rgb, 255:red, 0; green, 0; blue, 0 }  ][line width=0.75]    (0,5.59) -- (0,-5.59)   ;
\draw    (26.53,174.29) -- (214.26,174) ;
\draw [shift={(216.26,173.99)}, rotate = 539.9100000000001] [color={rgb, 255:red, 0; green, 0; blue, 0 }  ][line width=0.75]    (10.93,-3.29) .. controls (6.95,-1.4) and (3.31,-0.3) .. (0,0) .. controls (3.31,0.3) and (6.95,1.4) .. (10.93,3.29)   ;
\draw    (26.53,54.76) -- (26.53,114.52) ;
\draw [shift={(26.53,84.64)}, rotate = 270] [color={rgb, 255:red, 0; green, 0; blue, 0 }  ][line width=0.75]    (0,5.59) -- (0,-5.59)   ;
\draw    (26.53,114.52) -- (26.53,174.29) ;
\draw [shift={(26.53,144.41)}, rotate = 270] [color={rgb, 255:red, 0; green, 0; blue, 0 }  ][line width=0.75]    (0,5.59) -- (0,-5.59)   ;
\draw    (26.53,174.29) -- (86.95,174.23) ;
\draw [shift={(56.74,174.26)}, rotate = 539.94] [color={rgb, 255:red, 0; green, 0; blue, 0 }  ][line width=0.75]    (0,5.59) -- (0,-5.59)   ;
\draw    (56.74,174.26) -- (117.15,174.2) ;
\draw [shift={(86.95,174.23)}, rotate = 539.94] [color={rgb, 255:red, 0; green, 0; blue, 0 }  ][line width=0.75]    (0,5.59) -- (0,-5.59)   ;
\draw    (86.95,174.23) -- (147.36,174.17) ;
\draw [shift={(117.15,174.2)}, rotate = 539.94] [color={rgb, 255:red, 0; green, 0; blue, 0 }  ][line width=0.75]    (0,5.59) -- (0,-5.59)   ;
\draw    (117.15,174.2) -- (177.57,174.14) ;
\draw [shift={(147.36,174.17)}, rotate = 539.94] [color={rgb, 255:red, 0; green, 0; blue, 0 }  ][line width=0.75]    (0,5.59) -- (0,-5.59)   ;
\draw    (147.36,174.17) -- (207.77,174.11) ;
\draw [shift={(177.57,174.14)}, rotate = 539.94] [color={rgb, 255:red, 0; green, 0; blue, 0 }  ][line width=0.75]    (0,5.59) -- (0,-5.59)   ;
\draw [line width=2.25]  [dash pattern={on 2.53pt off 3.02pt}]  (61.5,33.39) -- (84.22,33.39) ;
\draw    (56.67,115.09) -- (56.74,174.26) ;
\draw   (54.07,115.09) .. controls (54.07,113.68) and (55.23,112.53) .. (56.67,112.53) .. controls (58.11,112.53) and (59.27,113.68) .. (59.27,115.09) .. controls (59.27,116.5) and (58.11,117.65) .. (56.67,117.65) .. controls (55.23,117.65) and (54.07,116.5) .. (54.07,115.09) -- cycle ;
\draw    (86.97,144.97) -- (86.95,174.23) ;
\draw   (84.37,144.97) .. controls (84.37,143.56) and (85.53,142.41) .. (86.97,142.41) .. controls (88.4,142.41) and (89.57,143.56) .. (89.57,144.97) .. controls (89.57,146.38) and (88.4,147.52) .. (86.97,147.52) .. controls (85.53,147.52) and (84.37,146.38) .. (84.37,144.97) -- cycle ;
\draw    (177.5,114.97) -- (177.57,174.14) ;
\draw [line width=2.25]  [dash pattern={on 2.53pt off 3.02pt}]  (177.52,88.27) -- (177.54,125.21) -- (177.57,174.14) ;
\draw   (174.92,85.71) .. controls (174.92,84.3) and (176.08,83.15) .. (177.52,83.15) .. controls (178.95,83.15) and (180.12,84.3) .. (180.12,85.71) .. controls (180.12,87.12) and (178.95,88.27) .. (177.52,88.27) .. controls (176.08,88.27) and (174.92,87.12) .. (174.92,85.71) -- cycle ;
\draw [line width=2.25]  [dash pattern={on 2.53pt off 3.02pt}]  (86.97,147.52) -- (86.95,174.23) ;
\draw [line width=2.25]  [dash pattern={on 2.53pt off 3.02pt}]  (56.67,117.65) -- (56.74,174.26) ;
\draw   (84.22,33.39) .. controls (84.22,31.98) and (85.39,30.84) .. (86.82,30.84) .. controls (88.26,30.84) and (89.42,31.98) .. (89.42,33.39) .. controls (89.42,34.8) and (88.26,35.95) .. (86.82,35.95) .. controls (85.39,35.95) and (84.22,34.8) .. (84.22,33.39) -- cycle ;
\draw    (61.3,48.78) -- (86.62,48.78) ;
\draw [shift={(86.62,48.78)}, rotate = 180] [color={rgb, 255:red, 0; green, 0; blue, 0 }  ][line width=0.75]    (0,5.59) -- (0,-5.59)   ;
\draw    (56.74,174.26) -- (56.67,115.09) ;
\draw [shift={(56.67,115.09)}, rotate = 449.93] [color={rgb, 255:red, 0; green, 0; blue, 0 }  ][line width=0.75]    (0,5.59) -- (0,-5.59)   ;
\draw    (86.95,174.23) -- (86.97,144.97) ;
\draw [shift={(86.97,144.97)}, rotate = 450.04] [color={rgb, 255:red, 0; green, 0; blue, 0 }  ][line width=0.75]    (0,5.59) -- (0,-5.59)   ;
\draw    (177.57,174.14) -- (177.5,114.97) ;
\draw [shift={(177.5,114.97)}, rotate = 449.93] [color={rgb, 255:red, 0; green, 0; blue, 0 }  ][line width=0.75]    (0,5.59) -- (0,-5.59)   ;
\draw    (300.62,173.46) -- (300.62,56.76) ;
\draw [shift={(300.62,54.76)}, rotate = 450] [color={rgb, 255:red, 0; green, 0; blue, 0 }  ][line width=0.75]    (10.93,-3.29) .. controls (6.95,-1.4) and (3.31,-0.3) .. (0,0) .. controls (3.31,0.3) and (6.95,1.4) .. (10.93,3.29)   ;
\draw [shift={(300.62,114.11)}, rotate = 450] [color={rgb, 255:red, 0; green, 0; blue, 0 }  ][line width=0.75]    (0,5.59) -- (0,-5.59)   ;
\draw    (300.62,173.46) -- (489.48,173.16) ;
\draw [shift={(491.48,173.16)}, rotate = 539.9100000000001] [color={rgb, 255:red, 0; green, 0; blue, 0 }  ][line width=0.75]    (10.93,-3.29) .. controls (6.95,-1.4) and (3.31,-0.3) .. (0,0) .. controls (3.31,0.3) and (6.95,1.4) .. (10.93,3.29)   ;
\draw    (300.62,54.76) -- (300.62,114.11) ;
\draw [shift={(300.62,84.43)}, rotate = 270] [color={rgb, 255:red, 0; green, 0; blue, 0 }  ][line width=0.75]    (0,5.59) -- (0,-5.59)   ;
\draw    (300.62,114.11) -- (300.62,173.46) ;
\draw [shift={(300.62,143.78)}, rotate = 270] [color={rgb, 255:red, 0; green, 0; blue, 0 }  ][line width=0.75]    (0,5.59) -- (0,-5.59)   ;
\draw    (300.62,173.46) -- (361.39,173.4) ;
\draw [shift={(331,173.43)}, rotate = 539.94] [color={rgb, 255:red, 0; green, 0; blue, 0 }  ][line width=0.75]    (0,5.59) -- (0,-5.59)   ;
\draw    (331,173.43) -- (391.78,173.37) ;
\draw [shift={(361.39,173.4)}, rotate = 539.94] [color={rgb, 255:red, 0; green, 0; blue, 0 }  ][line width=0.75]    (0,5.59) -- (0,-5.59)   ;
\draw    (361.39,173.4) -- (422.17,173.34) ;
\draw [shift={(391.78,173.37)}, rotate = 539.94] [color={rgb, 255:red, 0; green, 0; blue, 0 }  ][line width=0.75]    (0,5.59) -- (0,-5.59)   ;
\draw    (391.78,173.37) -- (452.56,173.31) ;
\draw [shift={(422.17,173.34)}, rotate = 539.94] [color={rgb, 255:red, 0; green, 0; blue, 0 }  ][line width=0.75]    (0,5.59) -- (0,-5.59)   ;
\draw    (422.17,173.34) -- (482.95,173.28) ;
\draw [shift={(452.56,173.31)}, rotate = 539.94] [color={rgb, 255:red, 0; green, 0; blue, 0 }  ][line width=0.75]    (0,5.59) -- (0,-5.59)   ;
\draw [line width=2.25]  [dash pattern={on 2.53pt off 3.02pt}]  (337.68,33.23) -- (360.54,33.23) ;
\draw    (330.93,114.67) -- (331,173.43) ;
\draw   (328.32,114.67) .. controls (328.32,113.27) and (329.49,112.13) .. (330.93,112.13) .. controls (332.38,112.13) and (333.55,113.27) .. (333.55,114.67) .. controls (333.55,116.07) and (332.38,117.21) .. (330.93,117.21) .. controls (329.49,117.21) and (328.32,116.07) .. (328.32,114.67) -- cycle ;
\draw   (358.8,144.34) .. controls (358.8,142.94) and (359.97,141.8) .. (361.41,141.8) .. controls (362.86,141.8) and (364.03,142.94) .. (364.03,144.34) .. controls (364.03,145.74) and (362.86,146.88) .. (361.41,146.88) .. controls (359.97,146.88) and (358.8,145.74) .. (358.8,144.34) -- cycle ;
\draw [line width=2.25]  [dash pattern={on 2.53pt off 3.02pt}]  (452.51,88.03) -- (452.53,124.72) -- (452.56,173.31) ;
\draw   (449.89,85.49) .. controls (449.89,84.09) and (451.06,82.96) .. (452.51,82.96) .. controls (453.95,82.96) and (455.12,84.09) .. (455.12,85.49) .. controls (455.12,86.9) and (453.95,88.03) .. (452.51,88.03) .. controls (451.06,88.03) and (449.89,86.9) .. (449.89,85.49) -- cycle ;
\draw [line width=2.25]  [dash pattern={on 2.53pt off 3.02pt}]  (361.41,146.88) -- (361.39,173.4) ;
\draw [line width=2.25]  [dash pattern={on 2.53pt off 3.02pt}]  (330.93,117.21) -- (331,173.43) ;
\draw   (360.54,33.23) .. controls (360.54,31.83) and (361.71,30.7) .. (363.15,30.7) .. controls (364.6,30.7) and (365.77,31.83) .. (365.77,33.23) .. controls (365.77,34.64) and (364.6,35.77) .. (363.15,35.77) .. controls (361.71,35.77) and (360.54,34.64) .. (360.54,33.23) -- cycle ;
\draw    (337.47,48.51) -- (362.95,48.51) ;
\draw [shift={(362.95,48.51)}, rotate = 180] [color={rgb, 255:red, 0; green, 0; blue, 0 }  ][line width=0.75]    (0,5.59) -- (0,-5.59)   ;
\draw    (331,173.43) -- (330.93,114.67) ;
\draw [shift={(330.93,114.67)}, rotate = 449.93] [color={rgb, 255:red, 0; green, 0; blue, 0 }  ][line width=0.75]    (0,5.59) -- (0,-5.59)   ;

\draw (47.77,180.6) node [anchor=north west][inner sep=0.75pt]  [font=\footnotesize]  {$X_{1}$};
\draw (77.56,180.6) node [anchor=north west][inner sep=0.75pt]  [font=\footnotesize]  {$X_{2}$};
\draw (107.36,180.6) node [anchor=north west][inner sep=0.75pt]  [font=\footnotesize]  {$X_{3}$};
\draw (137.65,180.6) node [anchor=north west][inner sep=0.75pt]  [font=\footnotesize]  {$X_{4}$};
\draw (167.95,181.1) node [anchor=north west][inner sep=0.75pt]  [font=\footnotesize]  {$X_{5}$};
\draw (99.73,23.44) node [anchor=north west][inner sep=0.75pt]    {$x$};
\draw (98.28,37.32) node [anchor=north west][inner sep=0.75pt]    {$\phi _{s_{1}( x)}( x)$};
\draw (4.65,135.43) node [anchor=north west][inner sep=0.75pt]    {$1$};
\draw (4.65,105.65) node [anchor=north west][inner sep=0.75pt]    {$2$};
\draw (4.65,76.37) node [anchor=north west][inner sep=0.75pt]    {$3$};
\draw (322.03,179.68) node [anchor=north west][inner sep=0.75pt]  [font=\footnotesize]  {$X_{1}$};
\draw (352.01,179.68) node [anchor=north west][inner sep=0.75pt]  [font=\footnotesize]  {$X_{2}$};
\draw (381.98,179.68) node [anchor=north west][inner sep=0.75pt]  [font=\footnotesize]  {$X_{3}$};
\draw (412.46,179.68) node [anchor=north west][inner sep=0.75pt]  [font=\footnotesize]  {$X_{4}$};
\draw (442.94,180.17) node [anchor=north west][inner sep=0.75pt]  [font=\footnotesize]  {$X_{5}$};
\draw (374.19,23.3) node [anchor=north west][inner sep=0.75pt]    {$x$};
\draw (372.83,37.07) node [anchor=north west][inner sep=0.75pt]    {$\phi _{s_{4}( x)}( x)$};
\draw (278.64,134.82) node [anchor=north west][inner sep=0.75pt]    {$1$};
\draw (278.64,105.25) node [anchor=north west][inner sep=0.75pt]    {$2$};
\draw (278.64,76.17) node [anchor=north west][inner sep=0.75pt]    {$3$};

\end{tikzpicture}
\caption{}
\end{figure}

Notice that the retractions $\phi_{s_k(x)}(x)$  $k$-times iteratively substract $1$ from the norm of the last nonzero coordinate. See the examples of figure 4, where the graph refers to the norms of the coordinates. 
It is clear that the sequence $(s_n(x))\subset \N\cup\{0\}$ is strictly decreasing until it reaches 0 and becomes constant. See \eqref{relk} of the next Lemma \ref{propseq} for a detailed analysis of this behaviour.

\begin{lemma}\label{propseq}
Let $x\in S$ and $k\in\N\cup\{0\}$. Then the following equations are satisfied:
\begin{enumerate}
\item $$\phi_{s_{k}(x)-1}(x)=\phi_{s_{k+1}(x)}(x),$$\label{equa}
\item $$s_j(\phi_{s_{k-j}(x)}(x))=s_k(x)\;\;\;\forall j\in\{0,\dots,k\},$$\label{rels}
\item$$\big( \varphi_s\circ\cdots\circ\varphi_{s_k(x)} \big)\big( \phi_{s_k(x)}(x) \big)=\varphi_{s_k(x)}\big(\phi_{s_k(x)}(x)\big)=\phi_{s_{k+1}(x)}(x)\;\;\;\forall s\ge s_{k+1}(x)+1,$$\label{comp}
\item$$s_{k+1}(x)=s_{k}(x)-q_{n(\phi_{s_{k}(x)}(x))}.$$\label{relk}
\end{enumerate}
\begin{proof}
By the definition of $s_{k+1}(x)$ we know that $\phi_{s_{k}(x)-1}(x)\in N_{s_{k+1}(x)}$. So
$$\begin{aligned}\phi_{s_{k}(x)-1}(x)&=\big(\varphi_{s_{k+1}(x)+1}\circ\cdots\circ\varphi_{s_{k}(x)}\big)\circ\phi_{s_{k}(x)-1}(x)\\&=\big(\varphi_{s_{k+1}(x)+1}\circ\cdots\circ\varphi_{s_{k}(x)}\big)\circ\big(\varphi_{s_{k}(x)}\circ\cdots\circ\varphi_{s(x)}\big)(x)=\phi_{s_{k+1}(x)}(x),\end{aligned}$$
which proves \eqref{equa}. To prove \eqref{rels} we proceed by induction in $k$.  If $k=0$ it is straightforward to see that $s_0(\phi_{s_0(x)}(x))=s_0(x)$. Inductive step from $k-1$ to $k$.  If $j\in\{1,\dots,k\}$ then we let $y=\phi_{s_{k-j}(x)}(x)$, so by the inductive assumption we know that
$$s_{j-1}(y)=s_{j-1}\big( \phi_{s_{k-1-(j-1)}(x)}(x) \big)=s_{k-1}(x).$$
In this case we compute
$$\begin{aligned}s_j(\phi_{s_{k-j}(x)}(x))&=s_j(y)=s(\phi_{s_{j-1}(y)-1}(y))=s(\phi_{s_{k-1}(x)-1}(\phi_{s_{k-j}(x)}(x)))\\&=s(\phi_{s_{k-1}(x)-1}(x))=s_k(x).\end{aligned}$$
Otherwise, if $j=0$ then by \eqref{equa}
$$\begin{aligned}s_0(\phi_{s_{k}(x)}(x))&=s(\phi_{s_{k-1}(x)-1}(x))=s_k(x),\end{aligned}$$
which finishes the proof of \eqref{rels}.

We turn our attention to the point \eqref{comp}. This  shows the amount of trivial local retractions that act in between the nontrivial ones.

Consider $s\ge s_{1}(x)+1$. By the definition of $s_1(x)$ we have  that $\varphi_{s(x)}(x)=\phi_{s(x)-1}(x)\in N_{s_1(x)}$.
Hence,
$$\big( \varphi_s\circ\cdots\circ\varphi_{s(x)} \big)(x)=\varphi_{s(x)}(x)=\big( \varphi_{s_1(x)+1}\circ\cdots\circ\varphi_{s(x)}(x)\big)=\phi_{s_1(x)}(x),$$
that is,
\begin{equation}\big( \varphi_s\circ\cdots\circ\varphi_{s(x)} \big)(x)=\varphi_{s(x)}(x)=\phi_{s_{1}(x)}(x)\;\;\;\forall s\ge s_{1}(x)+1\label{base}.\end{equation}
This is exactly the equation \eqref{comp} for $k=0$.
Now, for an arbitrary $k\in\N$, just take $y=\phi_{s_{k-1}(x)-1}(x)=\phi_{s_k(x)}(x)$ so that by definition $s_k(x)=s(y)$ and by \eqref{rels} it follows that $s_1(y)=s_{k+1}(x)$. Using \eqref{base} with $y$ instead of $x$ we finish the proof of $\eqref{comp}$.

Next, let us prove $\eqref{relk}$. If $y=\phi_{s_{k-1}(x)-1}(x)=\phi_{s_{k}(x)}(x)\neq0$ then using  property $\ref{iter})$ of Lemma \ref{proplemma} we have that
$$\begin{aligned}s_{k+1}(x)&=s\big( \phi_{s_k(x)-1}(x) \big)=s\big( \phi_{s_k(x)-1}(y) \big)=s\big( \phi_{s(y)-1}(y) \big)\\&=s(y)-q_{n(y)}=s_k(x)-q_{n(\phi_{s_{k}(x)}(x))}.\end{aligned}$$
In the case when $y=0$ \eqref{relk} is trivially true as $s_k(x)=s(y)=0$ and $q_0=0$.
\end{proof}
\end{lemma}

The main feature of the definition of $\phi_s: S\to S$ is that it is a discretization of the Lipschitz retraction defined in \cite{HM21} (Section 3). To be more precise, we need to recall some definitions of \cite{HM21}:

For every $m,s\in\N$  we define the function $f_m^s:X\to \R$ as
$$f_m^s(x)=r^s_m\left(1-\sum\limits_{i=1}^{m-1}\frac{||x_i||}{r^s_i}\right).$$
For a given $x\in \bigcup\limits_{n\in\N}P_n(X)$ and $s\in\N$ we set $m(x,s)=\max\Big\{k\in\{1,\dots,n(x)+1\}\;:\;\sum\limits_{i=1}^{k-1}\frac{||x_i||}{r^s_i}\le1\Big\}$. Next, we define the constant mapping $F^{0}:X\to\{0\}$, and we also define the map $F^s:X\to D^s$ as the unique extension of the map $\widetilde{F^s}:\bigcup\limits_{n\in\N}P_n(X)\to D^s$ given by
$$\widetilde{F^s}(x)=\begin{cases}P_{m(x,s)-1}(x)+\frac{x_{m(x,s)}}{||x_{m(x,s)}||}f^s_{m(x,s)}(x)\;\;\;&\text{if }m(x,s)\le n(x),\\
x&\text{if }m(x,s)=n(x)+1.\end{cases}$$
In \cite{HM21} it is proved that $\widetilde{F^s}$ is well-defined and Lipschitz, with a Lipschitz norm independent of $s$. We would like to point out that $F^s$ is defined in \cite{HM21} as a limit of retractions but the definition given here is equivalent (see page 13 of \cite{HM21}). Also, it is important to mention that
$$F^s\circ F^t=F^{min(s,t)}\;\;\;\forall s,t\in\N.$$

\begin{lemma}\label{lemmabaseseq}
For every $k\in\N$ and $x\in S$,
\begin{equation}\label{similk}||F^{s_k(x)}(x)-\phi_{s_k(x)}(x)||\le 6b.\end{equation}
\begin{proof}
If $x=0$ then the above is trivially true, so we may assume that $x\in S\setminus\{0\}$. We proceed by giving first an alternate definition of $\phi_{s_k(x)}$ and $F^{s_k(x)}$ which  will be later convenient for the proof of the Lemma. One has to view $k\in \N$ as the number of nontrivial local retractions that act on $x$ when it is mapped by $\phi_{s_k(x)}$. In fact $\sum\limits_{i=1}^{n(x)}||x_i||$ is the amount of nontrivial local retractions needed to map $x$ to the origin, see figure 4. We need a couple of definitions before we get on to the equivalent definitions of $\phi_{s_k(x)}$ and $F^{s_k(x)}$.\\
For every $x\in S\setminus\{0\}$ and $k\in\N$ such that $k\le \sum\limits_{i=1}^{n(x)}||x_i||$ we claim that there is only one pair $I,J\in\N\cup\{0\}$ satisfying that $1\le J\le||x_{n(x)-I}||$ and $k=J+\sum\limits_{i=1}^{I}||x_{n(x)-i+1}||$. We will use the notation $I=i(k,x)$ and $J=j(k,x)$, respectively. Indeed, if $I_1,J_1$ and $I_2,J_2$ were two suitable pairs for $x\in S\setminus\{0\}$ and $k\in\{1,\dots,\sum\limits_{i=1}^{n(x)}||x_i||\}$, then assuming $I_1+1\le I_2$ we obtain
$$\begin{aligned}k=&J_2+\sum\limits_{i=1}^{I_2}||x_{n(x)-i+1}||\ge J_2+\sum\limits_{i=1}^{I_1+1}||x_{n(x)-i+1}||\\\le& J_2+J_1+\sum\limits_{i=1}^{I_1}||x_{n(x)-i+1}||=k+J_1>k,\end{aligned}$$
which is a contradiction. So $I_1=I_2$ and therefore $J_1=J_2$, proving our claim.

\begin{figure}\label{IJ2}
\centering

\tikzset{every picture/.style={line width=0.75pt}} 

\begin{tikzpicture}[x=0.75pt,y=0.75pt,yscale=-1,xscale=1]

\draw    (26.53,174.29) -- (26.53,56.76) ;
\draw [shift={(26.53,54.76)}, rotate = 450] [color={rgb, 255:red, 0; green, 0; blue, 0 }  ][line width=0.75]    (10.93,-3.29) .. controls (6.95,-1.4) and (3.31,-0.3) .. (0,0) .. controls (3.31,0.3) and (6.95,1.4) .. (10.93,3.29)   ;
\draw [shift={(26.53,114.52)}, rotate = 450] [color={rgb, 255:red, 0; green, 0; blue, 0 }  ][line width=0.75]    (0,5.59) -- (0,-5.59)   ;
\draw    (26.53,174.29) -- (214.26,174) ;
\draw [shift={(216.26,173.99)}, rotate = 539.9100000000001] [color={rgb, 255:red, 0; green, 0; blue, 0 }  ][line width=0.75]    (10.93,-3.29) .. controls (6.95,-1.4) and (3.31,-0.3) .. (0,0) .. controls (3.31,0.3) and (6.95,1.4) .. (10.93,3.29)   ;
\draw    (26.53,54.76) -- (26.53,114.52) ;
\draw [shift={(26.53,84.64)}, rotate = 270] [color={rgb, 255:red, 0; green, 0; blue, 0 }  ][line width=0.75]    (0,5.59) -- (0,-5.59)   ;
\draw    (26.53,114.52) -- (26.53,174.29) ;
\draw [shift={(26.53,144.41)}, rotate = 270] [color={rgb, 255:red, 0; green, 0; blue, 0 }  ][line width=0.75]    (0,5.59) -- (0,-5.59)   ;
\draw    (26.53,174.29) -- (86.95,174.23) ;
\draw [shift={(56.74,174.26)}, rotate = 539.94] [color={rgb, 255:red, 0; green, 0; blue, 0 }  ][line width=0.75]    (0,5.59) -- (0,-5.59)   ;
\draw    (56.74,174.26) -- (117.15,174.2) ;
\draw [shift={(86.95,174.23)}, rotate = 539.94] [color={rgb, 255:red, 0; green, 0; blue, 0 }  ][line width=0.75]    (0,5.59) -- (0,-5.59)   ;
\draw    (86.95,174.23) -- (147.36,174.17) ;
\draw [shift={(117.15,174.2)}, rotate = 539.94] [color={rgb, 255:red, 0; green, 0; blue, 0 }  ][line width=0.75]    (0,5.59) -- (0,-5.59)   ;
\draw    (117.15,174.2) -- (177.57,174.14) ;
\draw [shift={(147.36,174.17)}, rotate = 539.94] [color={rgb, 255:red, 0; green, 0; blue, 0 }  ][line width=0.75]    (0,5.59) -- (0,-5.59)   ;
\draw    (147.36,174.17) -- (207.77,174.11) ;
\draw [shift={(177.57,174.14)}, rotate = 539.94] [color={rgb, 255:red, 0; green, 0; blue, 0 }  ][line width=0.75]    (0,5.59) -- (0,-5.59)   ;
\draw [line width=2.25]  [dash pattern={on 2.53pt off 3.02pt}]  (61.5,33.39) -- (84.22,33.39) ;
\draw    (56.67,115.09) -- (56.74,174.26) ;
\draw   (54.07,115.09) .. controls (54.07,113.68) and (55.23,112.53) .. (56.67,112.53) .. controls (58.11,112.53) and (59.27,113.68) .. (59.27,115.09) .. controls (59.27,116.5) and (58.11,117.65) .. (56.67,117.65) .. controls (55.23,117.65) and (54.07,116.5) .. (54.07,115.09) -- cycle ;
\draw    (86.97,144.97) -- (86.95,174.23) ;
\draw   (84.37,144.97) .. controls (84.37,143.56) and (85.53,142.41) .. (86.97,142.41) .. controls (88.4,142.41) and (89.57,143.56) .. (89.57,144.97) .. controls (89.57,146.38) and (88.4,147.52) .. (86.97,147.52) .. controls (85.53,147.52) and (84.37,146.38) .. (84.37,144.97) -- cycle ;
\draw    (177.5,114.97) -- (177.57,174.14) ;
\draw [line width=2.25]  [dash pattern={on 2.53pt off 3.02pt}]  (177.52,88.27) -- (177.54,125.21) -- (177.57,174.14) ;
\draw   (174.92,85.71) .. controls (174.92,84.3) and (176.08,83.15) .. (177.52,83.15) .. controls (178.95,83.15) and (180.12,84.3) .. (180.12,85.71) .. controls (180.12,87.12) and (178.95,88.27) .. (177.52,88.27) .. controls (176.08,88.27) and (174.92,87.12) .. (174.92,85.71) -- cycle ;
\draw [line width=2.25]  [dash pattern={on 2.53pt off 3.02pt}]  (86.97,147.52) -- (86.95,174.23) ;
\draw [line width=2.25]  [dash pattern={on 2.53pt off 3.02pt}]  (56.67,117.65) -- (56.74,174.26) ;
\draw   (84.22,33.39) .. controls (84.22,31.98) and (85.39,30.84) .. (86.82,30.84) .. controls (88.26,30.84) and (89.42,31.98) .. (89.42,33.39) .. controls (89.42,34.8) and (88.26,35.95) .. (86.82,35.95) .. controls (85.39,35.95) and (84.22,34.8) .. (84.22,33.39) -- cycle ;
\draw    (61.3,48.78) -- (86.62,48.78) ;
\draw [shift={(86.62,48.78)}, rotate = 180] [color={rgb, 255:red, 0; green, 0; blue, 0 }  ][line width=0.75]    (0,5.59) -- (0,-5.59)   ;
\draw    (56.74,174.26) -- (56.67,115.09) ;
\draw [shift={(56.67,115.09)}, rotate = 449.93] [color={rgb, 255:red, 0; green, 0; blue, 0 }  ][line width=0.75]    (0,5.59) -- (0,-5.59)   ;
\draw    (86.95,174.23) -- (86.97,144.97) ;
\draw [shift={(86.97,144.97)}, rotate = 450.04] [color={rgb, 255:red, 0; green, 0; blue, 0 }  ][line width=0.75]    (0,5.59) -- (0,-5.59)   ;
\draw    (177.57,174.14) -- (177.5,114.97) ;
\draw [shift={(177.5,114.97)}, rotate = 449.93] [color={rgb, 255:red, 0; green, 0; blue, 0 }  ][line width=0.75]    (0,5.59) -- (0,-5.59)   ;
\draw    (300.62,173.46) -- (300.62,56.76) ;
\draw [shift={(300.62,54.76)}, rotate = 450] [color={rgb, 255:red, 0; green, 0; blue, 0 }  ][line width=0.75]    (10.93,-3.29) .. controls (6.95,-1.4) and (3.31,-0.3) .. (0,0) .. controls (3.31,0.3) and (6.95,1.4) .. (10.93,3.29)   ;
\draw [shift={(300.62,114.11)}, rotate = 450] [color={rgb, 255:red, 0; green, 0; blue, 0 }  ][line width=0.75]    (0,5.59) -- (0,-5.59)   ;
\draw    (300.62,173.46) -- (489.48,173.16) ;
\draw [shift={(491.48,173.16)}, rotate = 539.9100000000001] [color={rgb, 255:red, 0; green, 0; blue, 0 }  ][line width=0.75]    (10.93,-3.29) .. controls (6.95,-1.4) and (3.31,-0.3) .. (0,0) .. controls (3.31,0.3) and (6.95,1.4) .. (10.93,3.29)   ;
\draw    (300.62,54.76) -- (300.62,114.11) ;
\draw [shift={(300.62,84.43)}, rotate = 270] [color={rgb, 255:red, 0; green, 0; blue, 0 }  ][line width=0.75]    (0,5.59) -- (0,-5.59)   ;
\draw    (300.62,114.11) -- (300.62,173.46) ;
\draw [shift={(300.62,143.78)}, rotate = 270] [color={rgb, 255:red, 0; green, 0; blue, 0 }  ][line width=0.75]    (0,5.59) -- (0,-5.59)   ;
\draw    (300.62,173.46) -- (361.39,173.4) ;
\draw [shift={(331,173.43)}, rotate = 539.94] [color={rgb, 255:red, 0; green, 0; blue, 0 }  ][line width=0.75]    (0,5.59) -- (0,-5.59)   ;
\draw    (331,173.43) -- (391.78,173.37) ;
\draw [shift={(361.39,173.4)}, rotate = 539.94] [color={rgb, 255:red, 0; green, 0; blue, 0 }  ][line width=0.75]    (0,5.59) -- (0,-5.59)   ;
\draw    (361.39,173.4) -- (422.17,173.34) ;
\draw [shift={(391.78,173.37)}, rotate = 539.94] [color={rgb, 255:red, 0; green, 0; blue, 0 }  ][line width=0.75]    (0,5.59) -- (0,-5.59)   ;
\draw    (391.78,173.37) -- (452.56,173.31) ;
\draw [shift={(422.17,173.34)}, rotate = 539.94] [color={rgb, 255:red, 0; green, 0; blue, 0 }  ][line width=0.75]    (0,5.59) -- (0,-5.59)   ;
\draw    (422.17,173.34) -- (482.95,173.28) ;
\draw [shift={(452.56,173.31)}, rotate = 539.94] [color={rgb, 255:red, 0; green, 0; blue, 0 }  ][line width=0.75]    (0,5.59) -- (0,-5.59)   ;
\draw [line width=2.25]  [dash pattern={on 2.53pt off 3.02pt}]  (337.68,33.23) -- (360.54,33.23) ;
\draw    (330.93,114.67) -- (331,173.43) ;
\draw   (328.32,114.67) .. controls (328.32,113.27) and (329.49,112.13) .. (330.93,112.13) .. controls (332.38,112.13) and (333.55,113.27) .. (333.55,114.67) .. controls (333.55,116.07) and (332.38,117.21) .. (330.93,117.21) .. controls (329.49,117.21) and (328.32,116.07) .. (328.32,114.67) -- cycle ;
\draw   (358.8,144.34) .. controls (358.8,142.94) and (359.97,141.8) .. (361.41,141.8) .. controls (362.86,141.8) and (364.03,142.94) .. (364.03,144.34) .. controls (364.03,145.74) and (362.86,146.88) .. (361.41,146.88) .. controls (359.97,146.88) and (358.8,145.74) .. (358.8,144.34) -- cycle ;
\draw [line width=2.25]  [dash pattern={on 2.53pt off 3.02pt}]  (452.51,88.03) -- (452.53,124.72) -- (452.56,173.31) ;
\draw   (449.89,85.49) .. controls (449.89,84.09) and (451.06,82.96) .. (452.51,82.96) .. controls (453.95,82.96) and (455.12,84.09) .. (455.12,85.49) .. controls (455.12,86.9) and (453.95,88.03) .. (452.51,88.03) .. controls (451.06,88.03) and (449.89,86.9) .. (449.89,85.49) -- cycle ;
\draw [line width=2.25]  [dash pattern={on 2.53pt off 3.02pt}]  (361.41,146.88) -- (361.39,173.4) ;
\draw [line width=2.25]  [dash pattern={on 2.53pt off 3.02pt}]  (330.93,117.21) -- (331,173.43) ;
\draw   (360.54,33.23) .. controls (360.54,31.83) and (361.71,30.7) .. (363.15,30.7) .. controls (364.6,30.7) and (365.77,31.83) .. (365.77,33.23) .. controls (365.77,34.64) and (364.6,35.77) .. (363.15,35.77) .. controls (361.71,35.77) and (360.54,34.64) .. (360.54,33.23) -- cycle ;
\draw    (337.47,48.51) -- (362.95,48.51) ;
\draw [shift={(362.95,48.51)}, rotate = 180] [color={rgb, 255:red, 0; green, 0; blue, 0 }  ][line width=0.75]    (0,5.59) -- (0,-5.59)   ;
\draw    (331,173.43) -- (330.93,114.67) ;
\draw [shift={(330.93,114.67)}, rotate = 449.93] [color={rgb, 255:red, 0; green, 0; blue, 0 }  ][line width=0.75]    (0,5.59) -- (0,-5.59)   ;
\draw   (167.72,86.36) .. controls (164.19,86.29) and (162.39,88.03) .. (162.32,91.56) -- (162.32,91.56) .. controls (162.23,96.61) and (160.41,99.11) .. (156.87,99.04) .. controls (160.41,99.11) and (162.13,101.67) .. (162.03,106.72)(162.08,104.44) -- (162.03,106.72) .. controls (161.97,110.25) and (163.71,112.05) .. (167.24,112.12) ;
\draw   (361.28,201.88) .. controls (361.27,206.55) and (363.6,208.88) .. (368.27,208.89) -- (397.13,208.92) .. controls (403.8,208.93) and (407.13,211.27) .. (407.12,215.94) .. controls (407.13,211.27) and (410.46,208.94) .. (417.13,208.95)(414.13,208.94) -- (445.99,208.98) .. controls (450.66,208.99) and (452.99,206.66) .. (453,201.99) ;
\draw   (369.2,170.4) .. controls (372.7,170.44) and (374.47,168.71) .. (374.5,165.21) -- (374.5,165.21) .. controls (374.56,160.22) and (376.34,157.74) .. (379.83,157.78) .. controls (376.34,157.74) and (374.62,155.22) .. (374.67,150.22)(374.64,152.47) -- (374.67,150.22) .. controls (374.71,146.73) and (372.98,144.96) .. (369.48,144.92) ;

\draw (47.77,180.6) node [anchor=north west][inner sep=0.75pt]  [font=\footnotesize]  {$X_{1}$};
\draw (77.56,180.6) node [anchor=north west][inner sep=0.75pt]  [font=\footnotesize]  {$X_{2}$};
\draw (107.36,180.6) node [anchor=north west][inner sep=0.75pt]  [font=\footnotesize]  {$X_{3}$};
\draw (137.65,180.6) node [anchor=north west][inner sep=0.75pt]  [font=\footnotesize]  {$X_{4}$};
\draw (167.95,181.1) node [anchor=north west][inner sep=0.75pt]  [font=\footnotesize]  {$X_{5}$};
\draw (99.73,23.44) node [anchor=north west][inner sep=0.75pt]    {$x$};
\draw (98.28,37.32) node [anchor=north west][inner sep=0.75pt]    {$\phi _{s_{1}( x)}( x)$};
\draw (4.65,135.43) node [anchor=north west][inner sep=0.75pt]    {$1$};
\draw (4.65,105.65) node [anchor=north west][inner sep=0.75pt]    {$2$};
\draw (4.65,76.37) node [anchor=north west][inner sep=0.75pt]    {$3$};
\draw (322.03,179.68) node [anchor=north west][inner sep=0.75pt]  [font=\footnotesize]  {$X_{1}$};
\draw (352.01,179.68) node [anchor=north west][inner sep=0.75pt]  [font=\footnotesize]  {$X_{2}$};
\draw (381.98,179.68) node [anchor=north west][inner sep=0.75pt]  [font=\footnotesize]  {$X_{3}$};
\draw (412.46,179.68) node [anchor=north west][inner sep=0.75pt]  [font=\footnotesize]  {$X_{4}$};
\draw (442.94,180.17) node [anchor=north west][inner sep=0.75pt]  [font=\footnotesize]  {$X_{5}$};
\draw (374.19,23.3) node [anchor=north west][inner sep=0.75pt]    {$x$};
\draw (372.83,37.07) node [anchor=north west][inner sep=0.75pt]    {$\phi _{s_{4}( x)}( x)$};
\draw (278.64,134.82) node [anchor=north west][inner sep=0.75pt]    {$1$};
\draw (278.64,105.25) node [anchor=north west][inner sep=0.75pt]    {$2$};
\draw (278.64,76.17) node [anchor=north west][inner sep=0.75pt]    {$3$};
\draw (111.2,90) node [anchor=north west][inner sep=0.75pt]    {$j( 1,x)$};
\draw (171.4,209) node [anchor=north west][inner sep=0.75pt]    {$i( 1,x) =0$};
\draw (369.6,221) node [anchor=north west][inner sep=0.75pt]    {$i( 4,x) =3$};
\draw (384.4,140.6) node [anchor=north west][inner sep=0.75pt]    {$j( 4,x)$};

\end{tikzpicture}
\caption{}
\end{figure}

The meaning  of $i(k,x)$ is the count of how many $x_i$'s are retracted to zero by $\phi_{s_k(x)}$, while $j(k,x)$ tells us the difference between the norm of $x$ and $\phi_{s_k(x)}(x)$ in their first different coordinate. See the examples shown in figure 5. \\
With these two new quantities we obtain a more convenient definition of both $\phi_{s_k(x)}$ and $F^{s_k(x)}$.

\begin{claim}\label{redef}For every $k\in\N$ and $x\in S$ such that $k\le\sum\limits_{i=1}^{n(x)}||x_i||$ the equalities
\begin{equation}\label{eqphi}\phi_{s_k(x)}(x)=P_{n(x)-i(k,x)-1}(x)+\Psi_{||x_{n(x)-i(k,x)}||-j(k,x)}(x_{n(x)-i(k,x)}),\end{equation}
\begin{equation}\label{eqF}F^{s_k(x)}(x)=P_{n(x)-i(k,x)-1}(x)+R_{||x_{n(x)-i(k,x)}||-j(k,x)}(x_{n(x)-i(k,x)}),\end{equation}
are satisfied.
\begin{proof}
Equation \eqref{eqphi} for k=1 follows inmediately after realizing that $i(1,x)=0$ and $j(1,x)=1$ for every $x\in S\setminus\{0\}$. Now we prove \eqref{eqF} for $k=1$. Since $s\big(\varphi_{s(x)}(x)\big)=s_1(x)$, using $\ref{sphere})$ of Lemma \ref{proplemma} with $\varphi_{s(x)}(x)$ instead of $x$ we obtain that
$$\sum\limits_{i=1}^{n(x)-1}\frac{||x_i||}{r_i^{s_1(x)}}+\frac{||x_{n(x)}||-1}{r_{n(x)}^{s_1(x)}}=1.$$
Then,
$$\sum\limits_{i=1}^{n(x)-1}\frac{||x_i||}{r_i^{s_1(x)}}\le \sum\limits_{i=1}^{n(x)-1}\frac{||x_i||}{r_i^{s_1(x)}}+\frac{||x_{n(x)}||-1}{r_{n(x)}^{s_1(x)}}=1.$$
Also, using now $\ref{sphere})$ of Lemma \ref{proplemma} with $x$ we get that
$$\sum\limits_{i=1}^{n(x)}\frac{||x_i||}{r_i^{s_1(x)}}>\sum\limits_{i=1}^{n(x)}\frac{||x_i||}{r_i^{s(x)}}=1.$$
This means that $m(x,s_1(x))=n(x)$. Now, as 
$$\sum\limits_{i=1}^{n(x)-1}\frac{||x_i||}{r_i^{s_1(x)}}+\frac{||x_{n(x)}||-1}{r_{n(x)}^{s_1(x)}}=1\;\;\Rightarrow\;\; \frac{||x_{n(x)}||-1}{r_{n(x)}^{s_1(x)}}=1-\sum\limits_{i=1}^{n(x)-1}\frac{||x_i||}{r_i^{s_1(x)}},$$
it follows that
$$\begin{aligned}F^{s_1(x)}(x)&=P_{n(x)-1}(x)+\frac{x_{n(x)}}{||x_{n(x)}||}f_{n(x)}^{s_1(x)}(x)\\&=P_{n(x)-1}(x)+\frac{x_{n(x)}}{||x_{n(x)}||}r_{n(x)}^{s_1(x)}\bigg(1-\sum\limits_{i=1}^{n(x)-1}\frac{||x_i||}{r_i^{s_1(x)}}\bigg)\\&=P_{n(x)-1}(x)+\frac{x_{n(x)}}{||x_{n(x)}||}r_{n(x)}^{s_1(x)}\bigg( \frac{||x_{n(x)}||-1}{r_{n(x)}^{s_1(x)}} \bigg)\\&=P_{n(x)-1}(x)+\frac{||x_{n(x)}||-1}{||x_{n(x)}||}x_{n(x)}=P_{n(x)-1}(x)+R_{||x_{n(x)}||-1}(x_{n(x)}).\end{aligned}$$
Hence, both \eqref{eqphi} and \eqref{eqF} are established for $k=1$. We proceed with the inductive step from $k-1$ to $k\ge2$. We focus on \eqref{eqphi} since the proof of \eqref{eqF} is essentially the same. Let $y=\phi_{s_{k-1}(x)}(x)$ so that $\phi_{s_k(x)}(x)=\phi_{s_1(y)}(y)$. If $y=0$ then both sides of the equation \eqref{eqphi} are $0$ so we may assume that $y\neq 0$. Then,
$$y=P_{n(x)-i(k-1,x)-1}(x)+\Psi_{||x_{n(x)-i(k-1,x)}||-j(k-1,x)}(x_{n(x)-i(k-1,x)}).$$
We distinguish two cases. First, let us assume that $||x_{n(x)-i(k-1,x)}||=j(k-1,x)$. In this case  $y=P_{n(x)-i(k-1,x)-1}(x)$ meaning that if $\exists i\in\N$ such that $n(y)+1\le i\le n(x)-i(k-1,x)-1$ then $||x_i||=0$. Equivalently, if $\exists i\in\N$ such that $i(k-1,x)+2\le i \le n(x)-n(y)$ then $||x_{n(x)-i+1}||=0$. With that in mind we deduce that
$$k-1=\sum\limits_{i=1}^{i(k-1,x)+1}||x_{n(x)-i+1}||\;\;\Rightarrow\;\; k=1+\sum\limits_{i=1}^{n(x)-n(y)}||x_{n(x)-i+1}||,$$
with $1\le n(y)=n(x)-(n(x)-n(y))$. It follows that $i(k,x)=n(x)-n(y)$ and $j(k,x)=1$ since $||x_{n(x)-(n(x)-n(y))}||=||x_{n(y)}||=||y_{n(y)}||\ge1$. Computing  $\phi_{s_k(x)}(x)$ we have that
$$\begin{aligned}\phi_{s_k(x)}(x)=&\phi_{s_1(y)}(y)=P_{n(y)-1}(y)+\Psi_{||y_{n(y)}||-1}(y_{n(y)})\\=&P_{n(x)-i(k,x)-1}(x)+\Psi_{||x_{n(x)-i(k,x)}||-j(k,x)}(x_{n(x)-i(k,x)}).\end{aligned}$$
In the second case  $1\le j(k-1,x)\le||x_{n(x)-i(k-1,x)}||-1$.  We claim that $i(k,x)=i(k-1,x)$ and $j(k,x)=j(k-1,x)+1$. Indeed,
$$k-1=j(k-1,x)+\sum\limits_{i=1}^{i(k-1,x)}||x_{n(x)-i+1}||\;\Rightarrow\;k=(j(k-1,x)+1)+\sum\limits_{i=1}^{i(k-1,x)}||x_{n(x)-i+1}||,$$
with $j(k-1,x)+1\le ||x_{n(x)-i(k-1,x)}||$. As $y=P_{n(x)-i(k-1,x)-1}(x)+\Psi_{||x_{n(x)-i(k-1,x)}||-j(k-1,x)}(x_{n(x)-i(k-1,x)})$ with $||x_{n(x)-i(k-1,x)}||-j(k-1,x)\ge1$ it follows that $n(y)=n(x)-i(k-1,x)=n(x)-i(k,x)$ and 
$$\begin{aligned}||y_{n(y)}||=&\big|\big| \Psi_{||x_{n(x)-i(k-1,x)}||-j(k-1,x)}(x_{n(x)-i(k-1,x)}) \big|\big|\\=&||x_{n(x)-i(k-1,x)}||-j(k-1,x)=||x_{n(x)-i(k,x)}||-j(k,x)+1.\end{aligned}$$
Now, it is only a matter of computing $\phi_{s_k(x)}(x)$ again,
$$\begin{aligned}\phi_{s_k(x)}(x)=&\phi_{s_1(y)}(y)=P_{n(y)-1}(y)+\Psi_{||y_{n(y)}||-1}(y_{n(y)})\\=&P_{n(x)-i(k,x)-1}(y)\\&+\Psi_{||x_{n(x)-i(k,x)}||-j(k,x)}(\Psi_{||x_{n(x)-i(k,x)}||-j(k,x)+1}(x_{n(x)-i(k,x)}))\\=&P_{n(x)-i(k,x)-1}(x)+\Psi_{||x_{n(x)-i(k,x)}||-j(k,x)}(x_{n(x)-i(k,x)}).\end{aligned}$$
\end{proof}
\end{claim}

We return to the proof of \eqref{similk}. In fact, taking into account Claim \ref{redef}, if $k>\sum\limits_{i=1}^{n(x)}||x_i||=t$ then $i(t,x)=n(x)-1$ and $j(t,x)=||x_{n(x)}||$. So $\phi_{s_t(x)}(x)=0=F^{s_t(x)}(x)$ and we are done as $\phi_{s_k(x)}(x)=\phi_{s_k(x)}(\phi_{s_t(x)}(x))=\phi_{s_k(x)}(0)=0=F^{s_k(x)}(F^{s_t(x)}(x))=F^{s_k(x)}(x)$. Otherwise, if $k\le\sum\limits_{i=1}^{n(x)}||x_i||$ then by Claim \ref{redef} and equality \eqref{similrad},
$$\begin{aligned}\big|\big|\phi_{s_k(x)}&(x)-F^{s_k(x)}(x)\big|\big|\\=&\big|\big|\big(\Psi_{||x_{n(x)-i(k,x)}||-j(k,x)}-R_{||x_{n(x)-i(k,x)}||-j(k,x)}\big)(x_{n(x)-i(k,x)})\big|\big|\le 6b.\end{aligned}$$
\end{proof}
\end{lemma}

\begin{remark}
If $X$ has a Schauder basis, then one may see from claim \ref{redef} that $\phi_{s_k(x)}(x)=F^{s_k(x)}(x)$ for every $k\in\N$.
\end{remark}

The following lemma is needed for the case when $s\neq s_k(x)$ in the statement of Lemma \ref{lemmabaseseq}.

\begin{lemma}\label{basic}
Suppose that $n\in\N$ and $t,s\in\N\cup\{0\}$ are such that $t\le s$. If $x-F^t(x)=\frac{x_n}{||x_n||}$ for some $x\in X$ with $||x_n||>0$, then there exists $\lambda\in[0,1]$ such that $\big( F^{s}-F^t \big)(x)=\lambda\frac{x_n}{||x_n||}$.
\begin{proof}
By density we may assume that $x\in \bigcup\limits_{n\in\N}P_n(X)\setminus\{0\}$, so there exists an $n(x)\in\N$. Let us assume first that $t\ge1$. Obviously, as $F^t(x)\neq x$ we deduce that $x\notin D^t$ so then $m(x,t)\le n(x)\in\N$. Now, looking at the definition of $F^t$ it is clear that
$$\frac{x_n}{||x_n||}=x-F^t(x)=\bigg(\sum\limits_{i=m(x,t)+1}^{n(x)}x_i\bigg)+x_{m(x,t)}\bigg(1-\frac{f^t_{m(x,t)}(x)}{||x_{m(x,t)}||}\bigg).$$
Since $m(x,t)\le n(x)$ by the definition of $m(x,t)$ it is deduced that $P_{m(x,t)}(x)\notin D^t$ so $1-\frac{f^t_{m(x,t)}(x)}{||x_{m(x,t)}||}>0$. Then, $m(x,t)=n(x)=n$ and
\begin{equation}\label{ineqtilde}x_{n}-\frac{x_{n}}{||x_{n}||}f^t_{n}(x)=\frac{x_n}{||x_n||},\end{equation}
by the linear independence of the blocks of the FDD. Clearly $m(x,s)\ge m(x,t)=n$. As $n(x)=n$ there are only two possibilities, namely $m(x,s)=n$ or $m(x,s)=n+1$, equivalently $x\in D^s$. In the case when $m(x,s)=n$ we have that
$$F^s(x)-F^t(x)=F^s_k(x)-F^t_k(x)=\frac{x_n}{||x_n||}\big( f_n^s(x)-f_n^t(x) \big).$$
Here just looking at the definitions and taking into account that $m(x,s)=n$ implies $\sum\limits_{i=1}^n\frac{||x_i||}{r_i^s}>1$, it is easy to see that $0\le f_n^t(x)\le f_n^s(x)<||x_n||$. Now by \eqref{ineqtilde} we have that $||x_n||-f_n^t(x)=1$.
This finishes the argument in this case since
$$0\le f_n^s(x)-f_n^t(x) =f_n^s(x)-||x_n||+1<1.$$
In the remaining case when $x\in D^s$, $F^s_k(x)=x$ so we have $F^s(x)-F^t(x)=x-F^t(x)=\frac{x_n}{||x_n||}$ and we are done.\\
Finally, if $t=0$ then $x-F^t(x)=x=\frac{x_n}{||x_n||}$. In particular, $||x_n|| =1$ and we may assume $s\ge1$ because $s=0$ is trivially true with $\lambda=0$. Again, $m(x,s)=n$ or $x\in D^s$. If $m(x,s)=n$ then $0\le f_n^s(x)<||x_n||=1$ and
$$F^s(x)-F^t(x)=F^s(x)=f_n^s(x)\frac{x_n}{||x_n||},$$
so the statement of the Lemma is satisfied. Otherwise, if $x\in D^s$ then again $F^s(x)-F^t(x)=x=\frac{x_n}{||x_n||}$.
\end{proof}
\end{lemma}

\begin{lemma}\label{simil1}
For every $s\in\N$,
$$||\phi_s(x)-F^s(x)||\le1+6b\;\;\;\forall x\in S.$$
\begin{proof}
Take $x\in S\setminus N_s$. It is clear that $s\in\{1,\dots,s_0(x)-1\}$ so there exists a $k\in\N\cup\{0\}$ such that $s\in\{s_{k+1}(x),\dots,s_k(x)-1\}$. Then, by Lemma \ref{propseq}, equation \eqref{comp},
$$\phi_s(x)=\varphi_{s+1}\circ\cdots\circ\varphi_{s(x)}(x)=\big( \varphi_{s+1}\circ\cdots\circ\varphi_{s_k(x)} \big)\big( \phi_{s_k(x)}(x) \big)=\phi_{s_{k+1}(x)}(x).$$
Thanks to equations \eqref{eqphi} and \eqref{eqF} we know that for every $x\in S$, $s\big(\phi_{s_k(x)}(x)\big)=s\big(F^{s_k(x)}(x)\big)$. Now, if we denote $y=F^{s_k(x)}(x)$ then, again by equation \eqref{eqF} we have that
$$\begin{aligned}y-F^{s_{1}(y)}(y)=&y-\big(P_{n(y)-1}(y)+R_{||y_{n(y)}||-1}(y_{n(y)})\big)\\=&y_{n(y)}-R_{||y_{n(y)}||-1}(y_{n(y)})=\frac{y_{n(y)}}{||y_{n(y)}||}\end{aligned}.$$
Then, thanks to Lemma \ref{basic} with $t=s_1(y)$ we know that there is a $\lambda\in[0,1]$ such that
$$F^s(x)-F^{s_{k+1}(x)}(x)=F^s(y)-F^{s_{1}(y)}(y)=\lambda \frac{y_{n(y)}}{||y_{n(y)}||}.$$
Finally,
$$\begin{aligned}||\phi_s(x)-F^s(x)||=&||\phi_{s_{k+1}(x)}(x)-F^s(x)||\\\le&||\phi_{s_{k+1}(x)}(x)-F^{s_{k+1}(x)}(x)||+||F^{s_{k+1}(x)}(x)-F^{s}(x)||\\\le& 6b+\lambda.\end{aligned}$$
\end{proof}
\end{lemma}

\begin{theorem}\label{theoFDD}
If $X$ is a Banach space with an FDD, then every spiderweb grid respect to the FDD has a retractional basis. 
\begin{proof}
First we claim that  for every $x\in X$ and $s\in \N$
\begin{equation}\label{claim1}||F^{s+1}(x)-F^{s}(x)||\le1.\end{equation}
In fact, it is enough to prove that $||x-F^s(x)||\le1$ for every $x\in D^{s+1}$ by the conmutativity of $F^s$. If $x\in D^s$ then the inequality is trivial since $x=F^s(x)$. Then, let us assume that $x\in D^{s+1}\setminus D^s$. As $x\in D^{s+1}$ we obtain the following,
\begin{equation}\label{firstFDD}\begin{aligned}\sum\limits_{i=m(x,s)}^{n(x)}||x_i||&\le r^{s+1}_{m(x,s)}\sum\limits_{i=m(x,s)}^{n(x)}\frac{||x_i||}{r_i^{s+1}}\\=&r_{m(x,s)}^{s+1}\bigg( \sum\limits_{i=1}^{n(x)}\frac{||x_i||}{r_i^{s+1}}-\sum\limits_{i=1}^{m(x,s)-1}\frac{||x_i||}{r_i^{s+1}} \bigg)\\\le& r_{m(x,s)}^{s+1}\bigg( 1-\sum\limits_{i=1}^{m(x,s)-1}\frac{||x_i||}{r_i^{s+1}} \bigg)=f_{m(x,s)}^{s+1}.\end{aligned}\end{equation}
Next, since $x\notin D^s$ we know that $1-\frac{f_{m(x,s)}^s(x)}{||x_{m(x,s)}||}>0$ so that we may compute using \eqref{firstFDD},
$$\begin{aligned} ||x-F^s(x)||=&\bigg|\bigg| \bigg(\sum\limits_{i=m(x,s)}^{n(x)}x_i\bigg)-\frac{x_{m(x,s)}}{||x_{m(x,s)}||}f_{m(x,s)}^s(x) \bigg|\bigg|\\\le&\bigg( \sum\limits_{i=m(x,s)+1}^{n(x)}||x_i|| \bigg)+\bigg|\bigg| x_{m(x,s)}\bigg( 1-\frac{f_{m(x,s)}^s(x)}{||x_{m(x,s)}||} \bigg) \bigg|\bigg|\\=&\bigg(\sum\limits_{i=m(x,s)}^{n(x)}||x_i||\bigg)-f_{m(x,s)}^s(x)\le f^{s+1}_{m(x,s)}(x)-f_{m(x,s)}^s(x)\\=&r_{m(x,s)}^{s+1}-r_{m(x,s)}^s. \end{aligned}$$
By \eqref{proprn} we know $r_{m(x,s)}^{s+1}=\frac{s+1}{s}r_{m(x,s)}^s$ so then 
$$||x-F^s(x)||\le r_{m(x,s)}^{s+1}-r_{m(x,s)}^s=\frac{r^s_{m(x,s)}}{s}=r^1_{m(x,s)}\le1$$
and \eqref{claim1} is proved. Now, we are going to proceed as in Corollary \ref{theospider}. Let us consider $S$ to be the $(a,b)$-spiderweb grid fixed at the beginning of the section. Denote $N_{-1}=\emptyset$ and $k_s=\#(N_s\setminus N_{s-1})$ for every $s\in\N\cup\{0\}$. Recall that
$$S=\{x_{(s,i)}\}_{\substack{n\in\N\cup\{0\}\\i=1,\dots,k_s}},$$
where $\{x_{(s,1)},\dots,x_{(s,k_s)}\}=N_{s}\setminus N_{s-1}$ is an arbitrary order of $N_s\setminus N_{s-1}$ for every $s\in\N\cup\{0\}$.
Order S lexicographically according to the previous notation and define the retractions of the retractional basis $\varphi_{(s,i)}:S\to S_{(s,i)}:=\{x_{(m,j)}\}_{(m,j)\le(s,i)}$ as
$$\varphi_{(s,i)}(x)=\begin{cases}\phi_s(x)\;\;&\text{if }\phi_s(x)\in S_{(s,i)},\\\phi_{s-1}(x)\;\;&\text{if }\phi_s(x)\notin S_{(s,i)}.\end{cases}$$
It is easy to see that $\varphi_{(s,i)}$ is a well defined retraction. Also, a straightforward computation gives $\varphi_{(s,i)}\circ\varphi_{(m,j)}=\varphi_{\min((s,i),(m,j))}$ so that it only remains to prove that $\varphi_{(s,i)}$ is Lipschitz. Indeed, by Lemma \ref{simil1} and equation \eqref{claim1}, if $x\in S$ then
$$||\phi_{s-1}(x)-F^s(x)||\le||\phi_{s-1}(x)-F^{s-1}(x)||+||F^{s-1}(x)-F^s(x)||\le6b+2$$
so it follows that
$$||\varphi_{(s,i)}(x)-F^{s}(x)||\le\max\{||\phi_s(x)-F^s(x)||,||\phi_{s-1}(x)-F^s(x)||\}\le 6b+2.$$
Then for every $x,y\in S$,
$$||\varphi_{(s,i)}(x)-\varphi_{(s,i)}(y)||\le 2(6b+2)+||F^s(x)-F^s(y)||\le \bigg(\frac{12b+4}{a}+||F^s||_{\text{Lip}}\bigg)||x-y||.$$
We are done since $||F^s||_{\text{Lip}}$ is independent of $s$ and $S$ was fixed but arbitrary.
\end{proof}
\end{theorem}

\begin{corollary}\label{corbasis}
Every grid respect to a Schauder basis of a Banach space has a retractional basis.
\end{corollary}

\begin{definition}
A Banach space $X$ is said to be griddable respect to a sequence of finite dimensional subspaces $(X_n)$ of $X$ whenever $(X_n)$ is an FDD of $X$ and for every sequence of $(a,b)$-nets $G_n\subset X_n$ the set $G=\bigcup\limits_{n\in\N}\sum\limits_{k=1}^nG_k$ is a net of $X$.

We will say that $X$ is griddable if it is griddable for some sequence $(X_n)$.
\end{definition}

If a Banach space $X$ has a $c_0$-like FDD  $(X_n)$, i.e. $X=\oplus_{c_0}\sum X_n$, then it is clearly griddable.
If $X$ has a Schauder basis and contains a copy of $c_0\hookrightarrow X$ then  it is griddable (with respect to possibly different Schauder basis) by Theorem 4.1 of \cite{DOS+08}.

\begin{corollary}\label{corFDD}
If $X$ is griddable, then every net $N$ of $X$ has a retractional basis. In particular, if $X$ is a Banach space with a Schauder basis containing a copy of $c_0$, or if $X$ has a $c_0$-like FDD $(X_n)$  then $\mathcal{F}(N)$ has  a retractional basis.
\end{corollary}

\begin{remark}
Let us point out that separable $L_1$-preduals, as well as subspaces of $c_0$ admitting an shrinking FDD, meet the conditions of Corollary \ref{corFDD}. In fact, by Theorem 1 of \cite{JZ72} every subspace of a quotient of $\big(\sum G_n\big)_{c_0}$ with shrinking FDD, where the $G_n$ are finite dimensional spaces, admits a $c_0$-like FDD.
\end{remark}

\bigskip

\printbibliography
\end{document}